%% file: MainBCRSentNov2010.tex
\documentclass{mcma_v02_modif}
\usepackage{array}
\usepackage{dsfont}
\title{A probabilistic algorithm approximating solutions of a singular PDE of porous media type}
\headlinetitle{A probabilistic algorithm approximating a singular
PDE  }

\authorone{Nadia Belaribi}
\addressone{ Laboratoire d'Analyse, G\'{e}om\'{e}trie et Applications (LAGA),
Universit\'{e} Paris 13,
 99, avenue Jean-Baptiste Cl\'{e}ment,
F-93430 Villetaneuse and
 ENSTA ParisTech, Unit\'{e} de Math\'{e}matiques
appliqu\'{e}es, 32, Boulevard Victor, F-75739 Paris Cedex 15 }
\countryone{France} \emailone{belaribi@math.univ-paris13.fr}

\authortwo{Fran\c{c}ois Cuvelier}
\addresstwo{Laboratoire d'Analyse, G\'{e}om\'{e}trie et Applications (LAGA),
 Universit\'{e} Paris 13, 99, avenue Jean-Baptiste Cl\'{e}ment
F-93430 Villetaneuse} \countrytwo{France}
\emailtwo{cuvelier@math.univ-paris13.fr}

\authorthree{Francesco Russo}
\addressthree{ ENSTA ParisTech, Unit\'{e} de Math\'{e}matiques
appliqu\'{e}es, 32, Boulevard Victor, F-75739 Paris Cedex 15,  INRIA
Rocquencourt and Cermics Ecole des Ponts et Chauss\'{e}es, Projet
MATHFI Domaine de Voluceau, BP 105 F-78153 Le Chesnay Cedex }
\countrythree{France}
\emailthree{francesco.russo@ensta-paristech.fr}

\abstract{The object of this paper is  a one-dimensional
 generalized porous media
equation (PDE) with possibly discontinuous coefficient $\beta$,
which is well-posed as an evolution problem in $L^1(\mathbb{R})$. In
some recent papers of Blanchard et alia and Barbu et alia,  the
solution was represented by the solution of a non-linear stochastic
differential equation in law if the initial condition is a bounded
integrable function. We first extend this result, at least when
$\beta$ is continuous and the initial condition is only integrable
with some supplementary technical assumption. The main purpose of
the article consists in introducing and implementing a stochastic
particle algorithm to approach the solution to (PDE) which also fits
in the case when $\beta$ is possibly irregular, to  predict some
long-time behavior of the solution and in comparing with some recent
numerical deterministic techniques. }



\acknowledgments{Part of the work was done during the stay of the
first and third named authors at the Bielefeld University (SFB 701
and BiBoS). They are grateful for the invitation.}

\begin{document}

%*******************************%
% Please input your paper here. %
% (e.g. \input xxx.tex)         %

\input Belaribi_CR.tex
%*******************************%

%*****************************************************%
% Please choose one of the following options:         %
%                                                     %
% VERSION A is for BIBTeX application.                %
% You have to add the names of your databases:        %
%                                                     %
%\bibliography{database}
%                                                     %
% IF you prefer VERSION B for the LaTeX standard      %
% bibliography environment, please use the same style %
% as produced by mcma.bst (see the examples in the    %
% instructions for authors).                          %
%                                                     %
%\begin{thebibliography}{99}
%  \bibitem{Argyros:2007} S. A. Argyros and P. Dodos,
%    \emph{Genericity and amalgamation of classes of  Banach spaces},
%    Adv. Math. 209 (2007), pp.~666--748.
%\end{thebibliography}
%                                                     %
%*****************************************************%

\end{document}

%% file: Belaribi_CR.tex
{\bf First Version: November 12th 2010}
\bigskip

\section{Introduction}

The main aim of this work is to construct and implement a
probabilistic algorithm which will allow us to approximate solutions
of a porous media type equation with monotone irregular coefficient.
Indeed, we are interested in the parabolic problem below:
\begin{equation}\label{EDP}
 \left\{
 \begin{array}{ccl}
  \partial_tu(t,x)&=&\frac{1}{2} \partial_{xx}^2 \beta\left(u(t,x)\right),\ \ t\in\left[0,+\infty\right[,\\
        u(0,x)&=&u_0(dx), \ \  \ x \in \mathbb{R},\\
\end{array}
\right.
\end{equation}
in the sense of distributions, where $u_0$ is an initial probability
measure. If $u_0$ has a density,  we will still denote it by the
same letter. We look for a solution of (\ref{EDP}) with time
evolution in $L^1(\mathbb{R})$. We formulate the following
assumption:

\textbf{Assumption(A)}

\begin{description}

    \item (i) $\beta:\mathbb{R}\rightarrow \mathbb{R}$ such that $\beta
    |_{\mathbb{R}_+}$ is monotone.

    \item (ii) $\beta(0)=0$ and $\beta$ continuous at zero.

    \item  (iii) We assume the existence of $\lambda>0$ such that $(\beta + \lambda
    id)(\mathbb{R}_+)=(\mathbb{R}_+)$, $id(x)\equiv x$.
\end{description}

A monotone function $\beta_0:\mathbb{R}\rightarrow \mathbb{R}$ can
be completed into a graph by setting
$\beta_0(x)=[\beta_0(x_{-}),\beta_0(x_{+})]$. An odd function
$\beta_0:\mathbb{R}\rightarrow \mathbb{R}$ such that $\beta
    |_{\mathbb{R}_+}=\beta_0
    |_{\mathbb{R}_+}$ produces in this way a maximal monotone
    graph.

In this introduction, however $\beta$ and $\beta_0$ will be
considered single-valued for the sake of simplicity. We leave more
precise formulations (as in Proposition \ref{EDS3.1} and Theorem
\ref{TheoNLSDE}) for the body of the article.

We remark that if $\beta$ fulfills  Assumption(A), then the odd
symmetrized $\beta_0$ fulfills the more natural

\textbf{Assumption(A')}

\begin{description}
    \item (i) $\beta_0:\mathbb{R}\rightarrow \mathbb{R}$  is monotone.

    \item (ii) $\beta_0(0)=0$ and $\beta_0$ continuous at zero.

    \item  (iii) We assume the existence of $\lambda>0$ such that
 $(\beta_0 + \lambda
    id)(\mathbb{R})=(\mathbb{R})$, $id(x)\equiv x$.
\end{description}
We define $\Phi:\mathbb{R}\rightarrow \mathbb{R}_+$, setting
\begin{equation}\label{PhiFunct}
\Phi(u)=\left\{
          \begin{array}{ll}
             \sqrt{\frac{\displaystyle{\beta_0(u)}}{\displaystyle{u}}}
 & \hbox{ if  $u \neq 0$,} \\ \\
            C & \hbox{  if  $u=0$,}
          \end{array}
        \right.
\end{equation}
where $C \in [\underset{u \rightarrow 0^+}{\liminf}
\ \Phi(u),\underset{u \rightarrow 0^+}{\limsup}\  \Phi(u)]$.\\

Note that when $\beta(u)=u.|u|^{m-1} ,~m>1$, the partial
differential equation (PDE) in (\ref{EDP})
 is nothing else but the classical porous media equation.
In this case $\Phi(u) = |u|^\frac{m-1}{2}$ and in particular $C =
0$.

Our main target is to analyze the case of an irregular coefficient
$\beta$. Indeed, we are particularly interested in the case when
 $\beta$ is continuous excepted for a
possible jump at one positive point, say $u_c>0$. A typical example
is:
\begin{eqnarray} \label{Heav}
\beta(u)=H(u-u_c).u,
\end{eqnarray}
$H$ being the Heaviside function and  $u_c$ will be  called  {\it
critical value} or {\it critical threshold}.

\begin{definition} \label{DefND}
\begin{description}
    \item i) We will say that the PDE in (\ref{EDP}), or $\beta$ is \emph{non-degenerate}
if there is a constant $c_0>0$ such that $\Phi \geq c_0$,  on each
compact of $\mathbb{R}_+$.

    \item ii) We will say that the PDE in (\ref{EDP}), or $\beta$ is
    \emph{degenerate} if $\lim\limits_{u \to 0^+}{\Phi(u)=0}$.
\end{description}

\end{definition}

\begin{remark}
\begin{description}
    \item i) We remark that $\beta$ is non-degenerate if and only if $\liminf\limits_{u \to 0^+}{\Phi(u)>0}$.

    \item ii) We observe that $\beta$ may be neither degenerate nor
    non-degenerate.
\end{description}
\end{remark}

Of course, $\beta$ in (\ref{Heav}) is degenerate. Equation
(\ref{Heav}) constitutes a model intervening in  some self-organized
criticality (often called SOC) phenomena, see \cite{1}  for a
significant monography on the subject. We mention the interesting
physical paper \cite{2}, which makes reference to a system whose
evolution is similar to the evolution of a "snow layer"
 under the influence of an "avalanche effect" which starts whenever
 the top of the layer is bigger than a critical value $u_c$.

We, in particular, refer to \cite{3} (resp. \cite{10}),
 which concentrates on the avalanche phase  and therefore
 investigates the problem
(\ref{EDP}) discussing existence, uniqueness and probabilistic
representation when $\beta$ is non-degenerate (resp. degenerate).
The authors had in mind  the singular PDE in (\ref{EDP})  as a
macroscopic model for which they gave a microscopic view via a
 probabilistic representation provided by a
non-linear stochastic differential equation (NLSDE); the stochastic
equation is supposed to describe the evolution of a single point of
the layer. The analytical assumptions formulated  by the  authors
were Assumption(A) and the Assumption(B) below which postulates
linear growth for $\beta$.

\textbf{Assumption(B)}

There exists a constant $c>0$ such that $|\beta(u)| \le c|u|$.
%\end{assumption}

Obviously we have,

%\begin{assumption} \label{AssumpB'}
\textbf{Assumption(B')}

There exists a constant $c>0$ such that $|\beta_0(u)| \le c|u|$.
%\end{assumption}
Clearly  \eqref{Heav} fulfills Assumption (B).

To the best of our knowledge the first author who considered a
probabilistic representation (of the type studied in this paper) for
the solutions of non linear deterministic partial differential
equations was McKean \cite{7}.
%particularly in relation with the so called propagation of chaos.
 However, in his case, the coefficients were  smooth. From then on the literature steadily grew and
 nowadays there is a vast amount of contributions to the subject.

A probabilistic interpretation of \eqref{EDP} when
$\beta(u)=u.|u|^{m-1},~m>1$  was provided in \cite{8}. For the same
$\beta$, though the method could be adapted to the case where
$\beta$ is Lipschitz, in \cite{9},  the author  studied the
evolution problem \eqref{EDP} when the initial condition and the
evolution takes values in the class of probability distribution
functions on $\mathbb{R}$. He studied both the probabilistic
representation and the so-called \emph{propagation of chaos}.

At the level of probabilistic representation, under Assumptions(A)
and (B), supposing that $u_0$ has a bounded density, \cite{3} (resp.
\cite{10}) proves existence and uniqueness (resp.  existence) in law
for (NLSDE). In the present work we are interested in some
theoretical complements, but the main purpose consists
 in  examining  numerical implementations
provided by (NLSDE), in comparison  with  numerical deterministic
schemes appearing in one recent paper, see \cite{17}.

 Let us now describe the principle of the probabilistic representation.
The stochastic differential equation (in law) rendering the
probabilistic representation is given by the following  (NLSDE):
\begin{equation}\label{EDS}
 \left\{
 \begin{array}{ccl}
 Y_t&=&Y_0+ \int\limits_0^t \Phi(u(s,Y_s))dW_s,\\
        u(t,\cdot)&=& \mbox{Law density of} \  Y_t ,\ \ \forall t > 0,\\
        u(0,\cdot)&=& \mbox{$u_0\ \ \ \ $    Law  of } Y_0,
\end{array}
\right.
\end{equation}
where $W$ is a classical Brownian motion. The solution of that
equation may be visualized as a continuous process $Y$ on some
filtered probability space $(\Omega,
\mathcal{F},({\mathcal{F}}_t)_{t \geq 0}, \mathbb{P})$ equipped with
an  $({\mathcal{F}}_t)_{t \geq 0}$-Brownian motion $W$.

Until now, theoretical results about well-posedness (resp.
existence) for (\ref{EDS}) were established when $\beta$ is
non-degenerate (resp. possibly degenerate) and in the case when $u_0
\in \left(L^1\bigcap L^{\infty}\right)(\mathbb{R})$. Even if the
present paper concentrates on numerical experiments, two theoretical
contributions are performed when $\Phi$ is continuous.

$\bullet$ Initially our aim was to produce an algorithm which allows
to start even with a measure or an unbounded function as intial
condition. Unfortunately, up to now, our implementation techniques
do not allow to treat this case.

 A first significant theoretical
contribution is Theorem \ref{EDS_Prop} which consists in fact in
extending the probabilistic representation obtained by \cite{10}
 to the case when $u_0 \in L^1(\mathbb{R})$, locally
of bounded variation outside a discrete set of points.

$\bullet$ A second contribution consists in showing in the
non-degenerate case that the mollified version of PDE in \eqref{EDP}
is
    in fact equivalent to its probabilistic representation, even
    when the initial condition $u_0$ is a probability measure.
This is done in Theorem \ref{TheoMollif}.

The connection between (\ref{EDS}) and (\ref{EDP}) is then given by
the following result.

\begin{proposition} \label{PI1.3}
Let us assume the existence of a solution $Y$ for (\ref{EDS}). Let
$u(t,\cdot)$ be the law density of $Y_t$, $t > 0$,
  that we suppose to exist.

Then $u:\left[0,T\right]\times \mathbb{R} \rightarrow \mathbb{R}_+$
provides a solution in the sense of distributions of (\ref{EDP})
with $u_0=u(0,\cdot)$.
\end{proposition}
The proof is well-known, but we recall here the basic argument for
illustration purposes.
\begin{proof}
Let $\varphi \in C_0^{\infty}(\mathbb{R})$, $Y$ be a solution of the
problem (\ref{EDS}). We apply It\^o's formula to $\varphi(Y)$ to
obtain :
\[\varphi(Y_t)=\varphi(Y_0) +\int_0^t \varphi'(Y_s)\Phi(u(s,Y_s))dW_s+\frac{\displaystyle{1}}{\displaystyle{2}}\int_0^t \varphi''(Y_s)\Phi^2(u(s,Y_s))ds.\]
Taking the expectation we get :
\[\int_{\mathbb{R}}\varphi(y)u(t,y)dy=\int_{\mathbb{R}}\varphi(y)u_0(y)dy +\frac{\displaystyle{1}}{\displaystyle{2}}\int_0^t ds
\int_{\mathbb{R}} \varphi''(y)\Phi^2(u(s,y))u(s,y)dy.\] Using then
integration by parts and the expression of $\beta$, the expected
result follows.
\end{proof}

In the literature there are several contributions about
approximation of non-linear PDE's of parabolic type using a
stochastic particles system, with study of the chaos propagation. We
recall that the chaos propagation takes place if the components of a
vector describing the interacting particle system become
asymptotically independent, when the number of particles goes to
infinity. Note that, physically motivated applications can  be
found, for instance in numerical studies in hydro- or
plasma-physics; \cite{Dawson} and \cite{Hockney_Livre} are
contributions expressing a heuristic or formal point of view.

When the non-linearity is of the first order, a significant
contribution was given by  \cite{12};
\cite{BossyTalay95,BossyTalay97} performed the rate of convergence,
\cite{13} provided a chaos propagation result. We also quote
\cite{Cald_Pulvi}, where authors provided a propagation of chaos
result for the Burger's equation.

 In the case of porous media type
equation in \eqref{EDP} with  $\beta$ Lipschitz, \cite{11}
investigated the probabilistic representation for \eqref{EDP} and a
mollified related equation. There, the authors provided a rigorous
proof of propagation of chaos in the case of Lipschitz coefficients,
see Proposition 2.3, Proposition 2.5 and Theorem 2.7 of \cite{11} .

Outside the Lipschitz case, an alternative method
 for studying convergence was investigated by
 \cite{Oelsch, Oelsh87, Oelsch_Sim}, whose limiting PDEs concerned  a class of equations
including the case $\beta(u) = u + u^2, \quad u \ge 0 $. In fact
\cite{Oelsch_Sim} computed the numerical solution of a viscous
porous medium equation through a particle algorithm and studied the
 $L^2$-convergence rate to the analytical solution.
More recent papers concerning the chaos propagation when  $\beta(u)
= u^2$ first and $\beta(u) = \vert u \vert^{m-1}
 u,  m > 1$ was
proposed in \cite{Philipow} and \cite{Philipow_Figal}.

As far as the coefficient $\beta$ is discontinuous, at our
knowledge, up to now,
 there are no such results.
As we announced, we are particularly interested in an empirical
investigation of the stochastic particle algorithm  approaching the
solution $u$ of \eqref{EDP}
  at some instant $t$,
in several situations with regular or irregular coefficient. We
recall that $u(t, \cdot)$ is a  probability density. That algorithm
involves Euler schemes of stochastic differential equations,
Monte-Carlo simulations expressing the empirical law and
non-parametric density estimation of $u(t, \cdot)$ using Gaussian
kernels, see \cite{14} for an introduction to the kernel method.
This technique crucially depends on the window width $\varepsilon$
of the smoothing kernel. Classical statistical tools for choosing
that parameter are described for instance in \cite{14}, where the
following formula for choosing the optimal bandwidth $\varepsilon$,
in the sense of minimizing the asymptotic  \emph{ mean integrated
squared error} (MISE), is given by
\begin{equation}\label{OptimalEpsilon}
\varepsilon_t={\left(2n\sqrt{\pi}\|\partial_{xx}^2
u(t,\cdot)\|^2\right)}^{-\frac{1}{5}},
\end{equation}
where,  $n$ is the sample size and $\|\cdot\|$ denotes the classical
$L^2(\mathbb{R})$ norm.

Of course, the above expression does not yield an immediately
practicable method for choosing the optimal $\varepsilon$ since
\eqref{OptimalEpsilon} depends on the second derivative of the
density  $u$, which we are trying indeed to estimate.
%which we are trying to estimate in the first place.
Therefore, several techniques were proposed to get through this
problem. First, a natural and easy approach, often called the
\emph{rule of thumb}, replaced the target density $u$ at time $t$
 in the functional $\|\partial_{xx}^2
u\|$, by a reference distribution function. For instance, \cite{14}
assumed that the unknown density is a standard normal function and
obtained the following practically used formula
\begin{equation} \label{SilvEpsilon}
 \varepsilon = \left(\frac{\displaystyle{4}}{\displaystyle{3n}}\right)^{\frac{1}{5}} \hat{\sigma},
\end{equation}
  $\hat{\sigma}$ being the empirical standard deviation. A
version which is more robust to outliers in the sample, consists in
replacing $\hat{\sigma}$ by a measure of spread of the variance
involving  the interquartile range.
 For instance, see \cite{14} for detailed
computations.

The \emph{oversmoothing} methods rely on the fact that there is a
simple upper bound for the MISE-optimal bandwidth. In fact,
\cite{Terrel_1990}, gave a lower bound for the functional
$\|\partial_{xx}^2 u\|$ and thus an upper bound for $\varepsilon$ in
\eqref{OptimalEpsilon}; it proposed to use this upper bound as an
optimal window width, see also \cite{Terrell_Scott_1995} for
histograms.

The two  methods above seem to work well for unimodal densities.
However, they lead to arbitrarily bad estimates of the bandwidth
$\varepsilon$, when for instance, the true density is far from being
 Gaussian, especially when it is a multimodal law.

The \emph{least squares cross validation} (LSCV) method aimed to
estimate the bandwidth that minimizes the integrated squared error
(ISE), based on a "leave-one-out" kernel density  estimator, see
\cite{Rudemo_1982, Bowman_1984}. The problem is that, for the same
target distribution, the estimated bandwidth through different
samples has a big variance, which produces instability.

The \emph{biased cross-validation} (BCV) approach,  introduced in
\cite{Scott_Terrell_1987} minimizes the score function obtained by
replacing the functional $\|\partial_{xx}^2 u\|$ in the formula of
the MISE by an estimator $\|\partial_{xx}^2 \hat{u}\|$, where
$\hat{u}$ is the kernel estimator of $u$. In fact,
\cite{Scott_Terrell_1987} proposed the use of the minimizer of that
score function as optimal bandwidth. This method seemed to be more
stable than the LSCV but still has large bias. The slow rate of
convergence of both the LSCV and BCV approaches encouraged
significant research on faster converging methods.

A popular approach, commonly called \emph{plug-in} method, makes use
of an indirect  estimator of the density functional
$\|\partial_{xx}^2 u\|$ in formula \eqref{OptimalEpsilon}. This
technique comes back to the early paper \cite{Woodroofe};
%The first paper that probably
% implemented this technique is \cite{Woodroofe},
% which proposed it for estimating the
%density at a given point.
in this framework the estimator of $\|\partial_{xx}^2 u\|$ requires
the computation of a \emph{pilot bandwidth} $h$, which is quite
different from  the window width $\varepsilon$ used for the kernel
density estimate. Indeed, this optimal bandwidth $h$  depends on
unknown density functionals involving partial derivatives greater
than $2$.
% $\|\partial_{xxx}^3 u\|$.
 Following an idea of
\cite{Wand_Jones_1995}, one could express $h$ iteratively through
higher order derivatives. In this spirit, the natural associated
problem consists in estimating for some positive integer $s$, the
quantity
 $\|\partial_{x^s}^s u\|$,
%for some
%positive integer $s$
in terms of $\|\partial_{x^{s+\ell}}^{s+\ell} u\|$ for some positive
integer $\ell$; an  \emph{$\ell$-stage direct plug-in} approach may
consist in replacing the norm $\|\partial_{x^{s+\ell}}^{s+\ell} u\|$
by the norm of the $s + \ell$ derivative of a Gaussian density. In
the present paper  we implement this idea with $s = \ell = 2$.
Important contributions to that topic were
 \cite{SJ} and  \cite{Jones_al_1996}
who improved the method via the so-called "solve-the-equation"
plug-in method. By this technique, the pilot bandwidth $h$ used to
estimate $\|\partial_{xx}^2 u\|$, is written as a function of the
kernel bandwidth $\varepsilon$. We shall describe in Section
\ref{Section:Num. Impl.} in details
 this bandwidth selection procedure
applied in the case of our probabilistic algorithm.

We point out, that a more recent tool was developed in
\cite{BotevPreparation} which improved the idea in
\cite{SJ,Jones_al_1996} in the sense that \cite{BotevPreparation}
did not postulate any normal reference rule. However, the numerical
experiments that we have performed using the Matlab routine
developed by the first author of \cite{BotevPreparation} have not
produced better results in the case when  $\beta$ is defined by
\eqref{Heav}.

In the paper we examine empirically the stochastic particle
algorithm for approaching the solution to the PDE in the case
$\beta(u) = u^3$ and in the case $\beta$ given by \eqref{Heav}. For
this more peculiar case, we compare the approximation with the one
obtained by one recent analytic deterministic numerical method.

Problems of the same type as  \eqref{EDP}, in the case when $\beta$
is Lipschitz but possibly degenerate, were extensively studied from
both the theoretical and numerical deterministic points of view. In
general, the numerical analysis of  \eqref{EDP} is difficult for at
least one reason: the appearance of singularities for compactly
supported solutions in the case of an irregular initial condition.
An usual technique  to approximate \eqref{EDP} involves implicit
discretization in time: it requires, at each time step, the
discretization of a nonlinear elliptic problem. However, when
dealing with nonlinear problems, one generally tries to linearize
them in order to take advantage of efficient linear solvers. Linear
approximation schemes based on the so-called non linear Chernoff's
formula with a suitable relaxation parameter and which arises in the
theory of nonlinear semi groups, were studied for instance in
\cite{Berger_Brezis_Rogers}. We also cite \cite{Kacur}, where the
authors approximated  degenerate parabolic problems including those
of porous media type. In fact, they used nonstandard
semi-discretization in time and applied a Newton-like iterations to
solve the corresponding elliptic problems. More recently, different
approaches based on kinetic schemes for degenerate parabolic systems
have been investigated in \cite{Aregba}. Finally a new scheme based
on the maximum principle and on the perturbation and regularization
approach was proposed in \cite{Pop_Yong}.

At the best of our knowledge, up to now, there are no analytical
methods dealing with the case when $\beta$ is given by \eqref{Heav}.
However, we are interested in a sophisticated approach developed in
\cite{17} and which appears to be best suited to describe the
evolution of singularities and efficient for computing discontinuous
solutions. In fact, \cite{17} focuses onto diffusive
\emph{relaxation schemes} for the numerical approximation of
nonlinear parabolic equations, see \cite{Jin,Jin_Levermore}, and
references therein. Those relaxed schemes are based on a suitable
semi-linear hyperbolic system with relaxation terms. Indeed, this
reduction is carried out in order to obtain  schemes that are easy
to implement. Moreover, with this approach it is possible to improve
such schemes by using different numerical approaches i.e. either
finite volumes, finite differences or  high order accuracy methods.

In particular, the authors in \cite{17} coupled ENO (Essentially Non
Oscillatory) interpolating algorithms for space discretization, see
\cite{Shu}, in order to deal with discontinuous solutions and
prevent the onset of spurious oscillation, with IMEX (implicit
explicit) Runge-Kutta schemes for time advancement, see
\cite{ParRusso}, to obtain a high order method. We point out that
\cite{17} studied convergence and stability of the scheme only in
the case when $\beta$ is Lipschitz but possibly degenerate and $u_0
\in L^1(\mathbb{R})$.

As a byproduct of numerical experiments we can forecast the longtime
behavior of $u(t,\cdot)$ where $(t,x) \mapsto u(t,x)$ is the
solution of the considered PDE. We can reasonably postulate that the
closure of  $\{ u \in L^1(\mathbb{R}), \ u \ge 0,\ \int_{\mathbb{R}}
u(x)dx =1
 \ | \ \beta(u) = 0 \}$ is a limiting set, provided it is not empty as in
the case $\beta(u) = u^3$.

The paper is organized as follows. Section \ref{Section: Exist. uni.
Results} is devoted to a survey of existence and uniqueness results
for both the deterministic
 problem  (\ref{EDP}) and the non-linear SDE (\ref{EDS}) rendering the probabilistic representation of (\ref{EDP}). We in particular, recall
  the results given by authors of \cite{3,10} and we provide some additional results in the case when the initial condition of (\ref{EDP})
    belongs to $L^1(\mathbb{R})$ but it is not necessarily bounded.

In Section \ref{Section: Some Complements}, we settle the
theoretical basis for the implementation of our probabilistic
algorithm. We first approximate the NLSDE (\ref{EDS}) by a mollified
version replacing $u(t,\cdot)$, the law density of $Y_t$, by a given
smooth function. We then construct an interacting particle system
for which we supposed that propagation of chaos result is verified.
We drive the attention on Theorem \ref{TheoMollif} which links the
mollified PDE \eqref{EDPM} with its probabilistic representation.

Section \ref{Section:Num. Impl.} is devoted to the numerical
procedure implementing the probabilistic algorithm. We first
introduce an Euler scheme to obtain a discretized version of the
interacting particles system defined in  Section \ref{Section: Some
Complements}. We then discuss the optimal choice of the window width
$\varepsilon$.

In Section \ref{Section:Det. App}, we describe the numerical
deterministic approach we use to simulate solutions of  \eqref{EDP}.
In fact, following \cite{17}, we first use finite differences and
ENO schemes for the space discretization, then we perform an
explicit Runge-Kutta scheme for time integration.

In Section \ref{Section: Numerical Experiments}, we proceed to the
validation of the algorithms. In fact, the first numerical
experiments discussed in that section concern the classical porous
media equation
 whose exact solution, in the case when the inial condition is a
 delta Dirac function,
 is explicitly given by the so-called \emph{Barrenblatt-Pattle} density, see \cite{16}. Then, we concentrate on the Heaviside case, i.e. with
  $\beta$ of the form (\ref{Heav}). In fact, we perform several test cases
 according  to the critical
 threshold $u_c$ and to the initial condition $u_0$. Finally, we conclude this section by  some
 considerations about the longtime behavior of
 solutions of  (\ref{EDP}).
\section{Existence and uniqueness results}\label{Section: Exist. uni.
Results} We start with some basic analytical framework. If
$f:\mathbb{R}\rightarrow \mathbb{R}$ is a bounded function we will
denote  $\|f\|_{\infty}=\sup\limits_{x\in \mathbb{R}}|f(x)|$. By
$S(\mathbb{R}) $ we denote the space of rapidly decreasing
infinitely differentiable functions $\varphi: \mathbb{R}\rightarrow
\mathbb{R}$. We denote by $\mathcal{M}(\mathbb{R})$ and
$\mathcal{M}_+(\mathbb{R})$ the set of finite measures and positive
finite  measures respectively.
\subsection{The deterministic PDE}
Based on some clarifications of some classical papers \cite{4, 5,
6}, \cite{3} states the following theorem about existence and
uniqueness in the sense of distributions (in a proper way).
\begin{proposition} \label{EDS3.1}
Let $u_0 \in \left(L^1\bigcap L^{\infty}\right)(\mathbb{R}), \quad
u_0 \ge 0$. We suppose the validity of  Assumptions (A) and (B).
 Then there is a unique solution in the sense of distributions $u \in (L^1\bigcap L^{\infty})(\left[0,T\right] \times \mathbb{R})$ of
\begin{equation}\label{Eq:EDS_3.1}
 \left\{
 \begin{array}{ccl}
  \partial_tu &\in & \frac{1}{2} \partial_{xx}^2 \beta(u),\\
        u(0,x)&=&u_0(x),\\
\end{array}
\right.
\end{equation}
in the sense that, there exists a unique couple $(u,\eta_u) \in
((L^1\bigcap L^{\infty}) (\left[0,T\right] \times \mathbb{R}))^2$
such that
\[\int u(t,x)\varphi(x)dx=\int u_0(x)\varphi(x)dx + \frac{\displaystyle{1}}{\displaystyle{2}} \int_0^t ds \int \eta_u(s,x) \varphi''(x)dx,
 \forall \varphi \in S(\mathbb{R})\]
and \[\eta_u(t,x) \in \beta(u(t,x))\ \ \mbox{for} \ dt \otimes
dx\mbox{-a.e.}\ \ (t,x) \in\left[0,t\right]\times \mathbb{R}\]
Furthermore, $||u(t,.)||_{\infty} \leq ||u_0||_{\infty}$ for every
$t \in \left[0,T\right]$ and there is a unique version of $u$ such
that $u \in C(\left[0,T\right];L^1(\mathbb{R}))~(\subset
L^1(\left[0,T\right]\times \mathbb{R}))$.
\end{proposition}
One significant difficulty of previous framework is that the
coefficient $\beta$ is discontinuous; this forces us to consider
$\beta$ as a multivalued function even though $u$ is single-valued.
Being $\beta$, in general, discontinuous it is difficult to imagine
the level of space regularity of the solution $u(t,\cdot)$ at time
$t$. In fact, Proposition 4.5 of \cite{10} says that almost surely
$\eta_u(t,\cdot)$ belongs $dt$-a.e in $H^1(\mathbb{R})$ if $u_0 \in
\left(L^1 \bigcap L^{\infty}\right)(\mathbb{R}) $. This helps in
some cases to visualize the behavior of $u(t,\cdot)$. The
proposition below makes some assertions when $\beta$ is of the type
of \eqref{Heav}, which constitutes our pattern situation.
\begin{proposition} \label{P22}

Let us suppose $u_0 \in \left(L^1 \bigcap
L^{\infty}\right)(\mathbb{R})$ and $\beta$ defined by \eqref{Heav}.
For $t \ge 0$, we denote by
\[E_t^0=\{x|\ u(t,x)=u_c\}.\]
 For almost all $t >0$,
\begin{enumerate}
    \item  $E_t^0$ has a non empty interior;

    \item every point of $E_t^0$ is either a local minimum or a
    local maximum.
\end{enumerate}

\end{proposition}

\begin{remark} \label{RP22}
The first point  of the previous proposition means that at almost
each time $t>0$, the function $u(t,\cdot)$ remains constant on some
interval.

The second point means that if the function $u(t,\cdot)$ crosses the
barrier $u_c$, it has first to stay constant for some time.
\end{remark}

\begin{proof}[Proof of Proposition~\ref{P22}]

For the sake of simplicity we fix $t>0$ such that $\eta_u(t,\cdot)
\in H^1(\mathbb{R})$ and we write $u=u(t,\cdot)$ ,
$\eta_u=\eta_u(t,\cdot)$.

\begin{enumerate}
    \item Since $\eta_u\in H^1(\mathbb{R})$ it is continuous, then the
    set $\mathcal{D}_0=\{x \in \mathbb{R}| \  \eta_u(x) \in
    ]0,u_c[\}$ is open. If $\eta_u(x)\in
    ]0,u_c[$ necessarily we have $u(x)=u_c$; in fact, if $u(x)<u_c$
    then $\eta_u(x)=0$ and if $u(x)>u_c$ then $\eta_u(x)=u(x)>u_c$.
    Since $\mathcal{D}_0$ is open and it is included in $E_t^0$ the
    result is established.

    \item Suppose the existence of sequences $(x_n)$ and $(y_n)$ such that $x_n \rightarrow x$
    with
    $u(x_n)<u_c$ and $y_n \rightarrow y$ with $u(y_n)>u_c$. By
    continuity of $\eta_u$ we have
    \[\eta_u(x_n)=0 {\underset{n \rightarrow \infty}{\rightarrow}}
    0=\eta_u(x)
\]
\[u(y_n)=\eta_u(y_n){\underset{n \rightarrow \infty}{\rightarrow}}
    \eta_u(x)=0,
\]
this is not possible because $u(y_n)>u_c$ for every $n$.
\end{enumerate}
\end{proof}
If $u_0 \in \mathcal{M}(\mathbb{R})$, we do not know any existence
or uniqueness theorem for (\ref{EDP}). Our first target consisted in
providing some generalization to Proposition \ref{EDS3.1} in the
case when $u_0$ is a finite measure. A solution in that case would
be, $u:\left]0,T\right] \times \mathbb{R} \rightarrow
L^1(\mathbb{R})$ continuous and such that
\[\lim_{t\rightarrow 0} u(t,dx)=u_0(dx),\]
weakly and where $u(t,dx)$ denotes $u(t,x)dx$. This is still an
object  of further technical investigations. For the moment, we are
only able to consider the case $u_0$ having a $ L^1(\mathbb{R})$
density still denoted by $u_0$, not necessarily bounded as in
Proposition \ref{EDS3.1}, at least when $\Phi$ characterized by
\eqref{PhiFunct} is continuous. In particular $\beta$ is also
continuous, but possibly degenerate. In that case, we can prove
existence of a
 distributional solution to (\ref{EDP}).
Even though this is not a very deep observation, this will settle
the basis of the corresponding probabilistic representation,
completely unknown in the literature. In fact, we provide the
following result.
\begin{proposition}\label{EDP_Prop}
Let $u_0 \in L^1(\mathbb{R})$. Furthermore, we suppose that
Assumption(A) and Assumption(B)  are fulfilled. We assume that
$\Phi$ is continuous on $\mathbb{R}_+$.

\begin{enumerate}
    \item [(1)] There is a solution $u$, in
the sense of distributions, to the problem
\begin{equation}\label{Eq1}
 \left\{
 \begin{array}{ccl}
  \partial_tu(t,x)&=&\frac{1}{2} \partial_{xx}^2 \beta(u(t,x)),~~t\in\left[0,\infty\right[,\\
        u(0,x)&=&u_0(dx),~~x \in \mathbb{R},\\
\end{array}
\right.
\end{equation}
in the sense that for every $\alpha \in S(\mathbb{R})$
\begin{equation}\label{Eq2}
\int_{\mathbb{R}}u(t,x)\alpha(x)dx=\int_{\mathbb{R}}u_0(x)\alpha(x)dx+
\frac{1}{2}\int_0^t ds\int_{\mathbb{R}}\alpha''(x)\beta(u(s,x))dx.
\end{equation}

    \item [(2)] If $u_0$ is locally of bounded variation excepted
    eventually on a discrete number of points $D_0$, then
    $\Phi(u(t,\cdot))$ has at most countable discontinuities  for every $t\in
    [0,T]$.
\end{enumerate}

\end{proposition}
\begin{proof}

 (1) Let $u_0 \in L^1(\mathbb{R})$, $u_0^N=u_0*\phi_{\frac{1}{N}},~N\in
\mathbb{N^\ast}$, where $\phi$ is a kernel with compact support and
$\phi_{\frac{1}{N}}(x) = N \phi(Nx), \ x \in  \mathbb{R} $. So
 $u^N_0$ is of class $C^1$, therefore
locally with bounded variation. Since $\|u_0^N\|_{\infty} \leq
\|\phi_{\frac{1}{N}}\|_{\infty}\|u_0\|_{L^1}$ then $u_0^N
\in(L^1\bigcap L^{\infty})(\mathbb{R})$. Moreover, we have
 \[\int_{\mathbb{R}}|u_0^N(x)-u_0(x)|dx \rightarrow 0,\  \mbox{as} \ \ N \rightarrow +\infty.\]
 On one hand, according to Proposition \ref{EDS3.1}, there is a unique solution
 $u^N$ of (\ref{Eq2}), i.e. for every
$\alpha \in S(\mathbb{R})$
\begin{equation}\label{etoile}
\int_{\mathbb{R}}u^N(t,x)\alpha(x)dx=\int_{\mathbb{R}}u_0^N(x)\alpha(x)dx+
\frac{1}{2}\int_0^t ds
\int_{\mathbb{R}}\alpha''(x)\beta(u^N(s,x))dx.
\end{equation}

On the other hand, according to Corollary 8.2 in Chap IV of
\cite{Show}, we have
\begin{equation}\label{E12}
\sup_{t\leq T}\int_{\mathbb{R}}|u^N(t,x)-u(t,x)|dx \rightarrow 0,\
\mbox{as} \ N\rightarrow +\infty.
\end{equation}
Therefore, there is a subsequence $(N_k)_{k\in\mathbb{N}}$ such that
\[u^{N_k}(t,x) \rightarrow u(t,x)\ \ dt\otimes dx\mbox{-a.e.},
\mbox{ as } \ k\rightarrow +\infty.  \] Since $\beta$ is continuous,
it follows that
\[\beta(u^{N_k}(t,x)) \rightarrow \beta(u(t,x))\ \ \ dt\otimes dx\mbox{-a.e.},
 \mbox{ as } \ k\rightarrow +\infty.  \]
Consequently, (\ref{etoile}) implies
\begin{equation}\label{etoile1}
\int_{\mathbb{R}}u(t,x)\alpha(x)dx=\int_{\mathbb{R}}u_0(x)\alpha(x)dx+
\lim_{k\rightarrow +\infty}\frac{1}{2}\int_0^t ds
\int_{\mathbb{R}}\alpha''(x)\beta(u^{N_k}(s,x))dx.
\end{equation}
In order to show that $u$ solves \eqref{Eq2}, we  verify
\begin{equation} \label{A11}
\lim_{N\to\infty} \int_0^t ds
\int_{\mathbb{R}}\alpha''(x)\beta(u^{N}(s,x))dx =\int_0^t ds
 \int_{\mathbb{R}}\alpha''(x)\beta(u(s,x))dx,
\end{equation}
where for notational simplicity we have replaced $N_k$ with $N$. So,
we can suppose that
\begin{equation} \label{A12}
u^N \rightarrow u, \quad \beta(u^N) \rightarrow \beta(u), \quad dt
\otimes dx\mbox{-a.e.}\quad \mbox{as} \ N\rightarrow +\infty.
 \end{equation}
Since $\vert \beta(u^N)\vert \le c \vert u^N \vert$ and $u^N
\rightarrow u $ in $L^1([0,T]\times \mathbb{R} )$, it follows that
$\beta(u^N)$ are equi-integrable. Consequently, by \eqref{A12},
 $\beta(u^N) \rightarrow \beta(u) $ in $L^1([0,T]\times \mathbb{R} )$,
and therefore \eqref{A11} follows. Finally, $u$ solves equation
(\ref{Eq2}).

(2) For this purpose we state a lemma concerning an elliptic
equation whose first statement item constitutes the kernel of the
proof of
 Proposition \ref{EDS3.1}.

Given $f:\mathbb{R}\rightarrow \mathbb{R}$, for $h \in \mathbb{R}$,
we denote \[f^h(x)=f(x+h)-f(x).\]
\begin{lemma}\label{Elliptic}
Let $f\in L^1$, $\lambda>0$.

\begin{enumerate}
\item There is a unique solution in the sense of distributions
of \[u-\lambda (\beta(u))''=f.\]

\item  Let $\chi$ be a smooth function with compact support. Then for
each $h$
\begin{eqnarray}\label{E3.9}
\int_{\mathbb{R}}\chi(x)|u^h(x)|dx \leq
\int_{\mathbb{R}}\chi(x)|f^h(x)|dx+ C\lambda|h|\|u\|_{L^1},
\end{eqnarray}
where $C$ is a constant depending on $\beta$ and $\chi$.
\end{enumerate}
\end{lemma}

\begin{proof}[Proof of Lemma~\ref{Elliptic}]

\emph{(i)} is stated in Theorem 4.1 of \cite{4} and Theorem 1 of
\cite{6}.

\emph{(ii)} The statement appears in Lemma 3.6 of \cite{10} in the
case when $f\in L^1\bigcap L^{\infty}$ but the proof remains the
same for $f \in L^1$.
\end{proof}

We go on with the proof of Proposition \ref{EDP_Prop}, point (2).
Let $\chi$ be a smooth nonnegative function with compact support on
$\mathbb{R}\backslash D_0$. We prove in fact
\begin{eqnarray}\label{LimSup}
\limsup_{h\to
0}\frac{\displaystyle{1}}{\displaystyle{h}}\int_{\mathbb{R}}\chi(x)|u^h(t,x)|dx
\leq \|u_0 \chi\|_{\mbox{\footnotesize{var}}}+
C\int_{[0,T]\times\mathbb{R}}|u(s,x)|dsdx,
\end{eqnarray}
where $\|\cdot\|_{\mbox{\footnotesize{var}}}$ denotes the total
variation and $C$ is a generic universal constant. For this purpose,
we proceed exactly as in the proof of Proposition 4.20 of \cite{10}
making use of Lemma \ref{Elliptic}. Inequality \eqref{LimSup}
allows, similarly as in \cite{10} to show that $u(t,\cdot)$
restricted to any compact interval of $\mathbb{R}\backslash D_0$ has
bounded variation. Therefore it has at most countable
discontinuities. Consequently, $\Phi(u(t,\cdot))$ has the same
property since $\Phi$ is supposed to be continuous.
\end{proof}
\subsection{The non-linear stochastic differential equation (NLSDE)}
\begin{definition}\label{DefExistProc}
We say that a process $Y$  is a solution to the NLSDE associated to
problem \eqref{EDP}, if there exists $\chi$ belonging to $
L^{\infty} (\left[0,T\right]\times \mathbb{R})$ such that;
\begin{equation}\label{NLSDE}
 \left\{
 \begin{array}{ccl}
  Y_t&=& Y_0+ \int_0^t \chi(s,Y_s)dW_s\ ,  \\
        \chi(t,x)\ &\in & \ \Phi(u(t,x)),\ \mbox{for}\ dt\otimes dx-\mbox{a.e.}\ (t,x)\in \left[0,T\right]\times \mathbb{R}, \\
u(t,x)\ &=&  \mbox{Law  density  of} \ Y_t, \ \ \ t > 0,\\
u(0,\cdot)\ &=& \ u_0,
\end{array}
\right.
\end{equation}
where $W$ is a Brownian motion  on some suitable filtered
probability space $(\Omega,\mathcal{F},(\mathcal{F}_t)_{t \geq
0},\mathbb{P})$. In particular, the first identity of \eqref{NLSDE}
holds in law. We introduce a notion appearing in \cite{10}.
\end{definition}

\begin{definition}\label{DefStrictIncr}
We say that $\beta$ is \emph{strictly increasing after some zero} if
there is a constant $c>0$ such that
\begin{description}
    \item i) $\beta|_{[0,c]}=0$.
    \item ii) $\beta$ is strictly increasing on $[c,+\infty[$.
\end{description}
\end{definition}
Up to now, two results are available concerning existence and
uniqueness of solutions to \eqref{NLSDE}. In fact, the first one is
stated in the case where $\beta$ is not degenerate and  the  second
one in the case when $\beta$ is degenerate, see respectively
\cite{3,10}. We summarize these two results in the following theorem
for easy reference later on.
\begin{theorem}\label{TheoNLSDE}
Let $u_0\in L^1\bigcap L^{\infty}$ such that $u_0 \geq 0$ and
$\int_{\mathbb{R}}u_0(x)dx=1$. Furthermore, we suppose that
Assumptions (A) and (B) are fulfilled.
\begin{enumerate}
    \item If $\beta$ is non-degenerate then it exists a solution $Y$
    to \eqref{NLSDE}, unique in law.

    \item Suppose $\beta$ is degenerate and either $\beta$ is strictly
    increasing after some zero or $u_0$ has locally bounded
    variation. Then there is a solution $Y$ not necessarily unique
    to \eqref{NLSDE}.
\end{enumerate}
\end{theorem}
A step forward is constituted by the proposition below. This
provides an existence result for the NLSDE, when $u_0$ is not
necessarily bounded at least whenever $\Phi$ is continuous. This
does not require a non-degenerate hypothesis on $\beta$.
\begin{theorem}\label{EDS_Prop}

Let $u_0 \in L^1(\mathbb{R})$  having  locally bounded variation
except on a discrete set of points $D_0$. Furthermore we suppose
that Assumption(A) and Assumption(B) are fulfilled. We assume that
$\Phi$ is continuous on $\mathbb{R}_+$.

The probabilistic representation related to \eqref{EDP} holds, i.e.
there is a process $Y$ solving \eqref{EDS} in law.
\end{theorem}
\begin{proof}

Let $u_0^N$ be the function considered at the beginning of the proof
of Proposition \ref{EDP_Prop}. According to Theorem \ref{TheoNLSDE},
let $Y_0^N$ be the solution to
\begin{equation}\label{EDPM1}
 \left\{
 \begin{array}{ccl}
 Y_t^N&=& Y_0^N+ \int_0^t \Phi(u^N(s,Y_s^N))dW_s,\\
        u^N(t,\cdot)&=& \mbox{Law density of } \ Y_t^N,\ \ \forall \ t \geq 0,\\
        u^N(0,\cdot)&=& u_0^N.
\end{array}
\right.
\end{equation}
Since $\Phi$ is bounded, using Burkholder-Davies-Gundy inequality
one obtains
\[E\left(Y_t^N-Y_s^N\right)^4\leq \mbox{const}(t-s)^2.\]
This implies ( see for instance Problem 4.11, Section 2.4 of
\cite{KARSH}) that the laws of $Y^N, N \geq 1$ are tight.
Consequently, there is a subsequence $Y^k:=Y^{N_k}$ converging in
law (as $C(\left[0,T\right])$-valued random elements) to some
process $Y$. We set $u^k=u^{N_k}$, where we recall that
$u^k(t,\cdot)$ is the law of $Y_t^k$. We also set
$X_t^k=Y_t^k-Y_0^k$. Since $[X^k]_t= \int\limits_0^t
\Phi^2(u^k(s,Y_s^k))ds$ and  $\Phi$ is bounded,
 the continuous local martingales $X^k$ are indeed martingales.

By Skorohod's theorem there is a new probability space
$(\widetilde{\Omega}, \widetilde{\mathcal{F}}, \widetilde{P})$ and
processes ${\widetilde{Y}}^k$, with the same distribution as $Y^k$
so that ${\widetilde{Y}}^k$ converges $\widetilde{P}$-a.s. to
 some process $\widetilde{Y}$, of course distributed as $Y$, as $C(\left[0,T\right])$-valued random element.
 In particular, the processes $\widetilde{X}_t^k=\widetilde{Y}_t^k-\widetilde{Y}_0^k$ remain martingales with respect to
  the filtration generated by  ${\widetilde{Y}}^k$. We denote the sequence ${\widetilde{Y}}^k$ (resp. $\widetilde{Y}$),
again by $Y^k$ (resp. Y).

We now aim to prove that
\begin{equation}\label{Eq5}
Y_t=Y_0+ \int_0^t \Phi(u(s,Y_s))dW_s,
\end{equation}
for some standard Brownian motion $W$ with respect with some
filtration $(\mathcal{F}_t)$.

We consider the stochastic process $X$ (vanishing at zero) defined
by $X_t=Y_t-Y_0$. We also set again $X_t^k=Y_t^k-Y_0^k$. Taking into
account Theorem 4.2 in Chap 3 of \cite{KARSH}, to establish
(\ref{Eq5}), it will be enough to prove that $X$ is an
$\mathcal{Y}$-martingale with quadratic variation $[X]_t=\int_0^t
\Phi^2(u(s,Y_s))ds$, where $\mathcal{Y}$ is the canonical filtration
associated with $Y$.

Let $s,t \in[0,T]$, with $t>s$ and $\psi$ a bounded continuous
function from $C([0,s])$ to $\mathbb{R}$. In order to prove the
martingale property for $X$, we need to show  that
\begin{equation}\label{YMart}
E\left[(X_t-X_s)\psi(Y_r,r \leq s)\right]=0.
\end{equation}
Since $Y^k$ are martingales, we have
\begin{equation}\label{YMartk}
E\left[(X_t^k-X_s^k)\psi(Y^k_r,r \leq s)\right]=0.
\end{equation}
Consequently \eqref{YMart} follows from \eqref{YMartk} and the fact
that $Y^k \rightarrow Y$  a.s. ($X^k \rightarrow X$ a.s.) as
$C(\left[0,T\right])$-valued random process. In fact for each $t
\geq 0,~X^k_t \rightarrow X_t$ in $L^1(\Omega)$ since $(X_t^k,k\in
\mathbb{N})$ is bounded in $L^2(\Omega)$.

It remains to show that $X_t^2-\int_0^t \Phi^2(u(s,Y_s))ds,\ \ \
t\in[0,T]$, defines an $\mathcal{Y}$-martingale, that is,  we need
to verify
\[E\left[\left(X_t^2-X_s^2-\int_s^t \Phi^2(u(r,Y_r))dr\right)\psi(Y_r,r \leq s)\right]=0.\]
We proceed similarly as in the proof of Theorem 4.3 in \cite{3} but
even with some simplification. For the comfort of the reader we give
a complete proof.

 The left-hand side decomposes into $I^1(k)+I^2(k)+I^3(k)$,
where
\begin{eqnarray*}
I^1(k)&=& E\left[\left(X_t^2-X_s^2-\int_s^t \Phi^2(u(r,Y_r))dr\right)\psi(Y_r,r \leq s)\right]\\
      &-& E\left[\left((X_t^k)^2-(X_s^k)^2-\int_s^t \Phi^2(u(r,Y_r^k))dr\right)\psi(Y_r^k,r \leq s)\right],\\\\
I^2(k)&=& E\left[\left((X_t^k)^2-(X_s^k)^2-\int_s^t \Phi^2(u^k(r,Y_r^k))dr\right)\psi(Y_r^k,r \leq s)\right],\\\\
I^3(k)&=& E\left[\left(\int_s^t (\Phi^2(u^k(r,Y_r^k))-\Phi^2(u(r,Y_r^k)))dr\right)\psi(Y_r^k,r \leq s)\right].\\
\end{eqnarray*}
 We start by showing the convergence of $I^3(k)$. Now, $\psi(Y_r^k,r \leq s)$ is dominated by a constant
 $C$. Clearly we have
 \begin{eqnarray*}
 I^3(k)\leq C\int_s^t
 dr\int_{\mathbb{R}}|\Phi^2(u^k(r,y))-\Phi^2(u(r,y))|u^k(r,y)dy.
 \end{eqnarray*}
 The right hand side of this inequality is equal to
 $C[J^1(k)+J^2(k)]$, where
 \begin{eqnarray*}
 J^1(k)&=&\int_s^t
 dr\int_{\mathbb{R}}|\Phi^2(u^k(r,y))-\Phi^2(u(r,y))|\left(u^k(r,y)-u(r,y)\right)dy,\\
 J^2(k)&=&\int_s^t
 dr\int_{\mathbb{R}}|\Phi^2(u^k(r,y))-\Phi^2(u(r,y))|u(r,y)dy.
 \end{eqnarray*}
 Since $u^k \rightarrow u$ in $C([0,T];L^1)$ and $\Phi^2$ is
 bounded then $\lim\limits_{k\to {+\infty}}J^1(k)=0$.

 Furthermore, there is a subsequence $(k_n)_{n\in \mathbb{N}}$ such
 that \[u^{k_n}(r,y) \rightarrow u(r,y) \ \ \ dr\otimes dy-\mbox{a.e.   as }  \ n \rightarrow +\infty.\]
 Since $\Phi^2$ is continuous, it follows that \[\Phi^2\left(u^{k_n}(r,y)\right) \rightarrow \Phi^2\left(u(r,y)\right) \ \ \ dr\otimes dy-\mbox{a.e.,  } N\rightarrow +\infty\]
 On the other hand, since \[|\Phi^2\left(u^{k_n}(r,y)\right)-\Phi^2\left(u(r,y)\right)|\leq 2\sup\limits_{u\in\mathbb{R}}{\Phi^2(u)}|u(r,y)|,\]
 Lebesgue's dominated convergence Theorem implies that $\lim\limits_{k\to
 {+\infty}}J^2(k)=0$.

Now we go on with the analysis of $I^2(k)$ and $I^1(k)$. $I^2(k)$
equals zero since $X^k$ is a martingale with quadratic variation
 given by $[X]_t=\int\limits_0^t\Phi^2(u^k(r,Y_r^k))dr$.

Finally, we treat $I^1(k)$. We recall that $X^k\rightarrow X$ a.s.
as a random element in $C([0,T])$ and that the sequence
$E((X_t^k)^4)$ is bounded, so $(X_t^k)^2$ are uniformly integrable.

Therefore, we have
\[E\left[\left((X_t)^2-(X_s)^2\right)\psi(Y_r,r \leq s)\right]-E\left[\left((X_t^k)^2-(X_s^k)^2\right)\psi(Y_r^k,r \leq s)\right]\rightarrow 0\]
when $k$ goes to infinity. It remains to prove that
\begin{equation}\label{Eq7}
\int_s^t
E\left[\left(\Phi^2(u(r,Y_r))-\Phi^2(u(r,Y_r^k))\right)\psi\left(Y_r,r
\leq s\right)dr\right]\rightarrow 0.
\end{equation}
Now, for fixed $dr$-a.e., $r\in[0,T]$, the set $\mathcal{S}(r)$ of
 discontinuities of  $\Phi(u(r,.))$
is countable   because of Proposition \ref{EDP_Prop}, point (2). The
law of $Y_r$ has a density and it is therefore non-atomic. Let
$N(r)$ be the  event of all $\omega \in \Omega$ such that
$Y_r(\omega)$ belongs to $\mathcal{S}(r).$ The probability of $N(r)$
equals $E(1_{ \mathcal{S}(r)} (Y_r)) = \int_\mathbb{R} 1_{
\mathcal{S}(r)}(y) dv(y)=0$, where $v$ is the law of  $Y_r$.
Consequently $N(r)$ is a negligible set.

For $\omega \notin N(r) $, we have $\lim\limits_{k \rightarrow
+\infty}
\Phi^2\left(u(r,Y_r^k(\omega))\right)=\Phi^2\left(u(r,Y_r(\omega))\right)$.
Since $\Phi$ is bounded, Lebesgue's dominated convergence theorem
implies (\ref{Eq7}).

Concerning the question wether $u(t,.)$ is the law of $Y_t$, we
recall that for all t,
 $Y_t^k$ converges (even in probability) to $Y_t$ and $u^k(t,.)$,
which is the law density of $Y_t^k$, goes to $u(t,.)$ in
$L^1(\mathbb{R})$. By the uniqueness of the limit in (\ref{Eq2}),
this obviously implies that $u(t,.)$ is the law density of $Y_t$.
\end{proof}
\section{Some complements related to the NLSDE}\label{Section: Some
Complements}
\subsection{A mollified version}

We suppose here that $u_0(dx)$ is a probability measure. Let $Y_0$
be a random variable distributed according to $u_0(dx)$ and
independent of the Brownian motion $W$.

In preparation to  numerical probability simulations, we define
$K_{\varepsilon }$ for every $\varepsilon>0$, as a smooth
regularization kernel obtained from a fixed probability density
function $K$ by the scaling :
\begin{equation}\label{Scaling}
K_{\varepsilon
}(x)=\frac{\displaystyle{1}}{\displaystyle{\varepsilon}}K\left(\frac{\displaystyle{x}}{\displaystyle{\varepsilon}}\right),\
\ x \in \mathbb{R}.
\end{equation}
We suppose in this section that $\Phi$ is single valued, therefore
continuous. This hypothesis will not be in force in Sections
\ref{Section:Num. Impl.} and \ref{Section: Numerical Experiments}.

In this subsection we wish to comment about the mollified version of
the NLSDE (\ref{EDS}), given by
\begin{equation}\label{EDSM}
 \left\{
 \begin{array}{ccl}
   Y_t^{\varepsilon}&=& Y_0+ \int\limits_0^t \Phi\left((K_{\varepsilon }*v^{\varepsilon})(s,Y_s^{\varepsilon})\right)dW_s,\\
        v^{\varepsilon}(t,\cdot)&=& \mbox{Law  of}
\ Y_t^{\varepsilon},\ \ \forall~t>0,\\
       v^{\varepsilon}(0,\cdot)&=& u_0
 \end{array}
\right.
\end{equation}
and its relation to the nonlinear integro-differential PDE
\begin{equation}\label{EDPM}
 \left\{
 \begin{array}{ccl}
\partial_tv^{\varepsilon}(t,x)&=&\frac{1}{2} \partial_{xx}^2 \left(\Phi^2(K_{\varepsilon}* v^{\varepsilon}(t,x)) v^{\varepsilon}(t,x)\right),\ (t,x)\in \left]0,+\infty\right[\times\mathbb{R},\\
        v^{\varepsilon}(0,\cdot)&=&u_0.
\end{array}
\right.
\end{equation}
 where, $t\mapsto v^{\varepsilon}(t,\cdot)$ may be measure-valued.

\begin{remark}\label{Remark_Mollif}
\begin{enumerate}
\item When $\Phi$ is Lipschitz, the authors of \cite{11} proved in
Proposition 2.2, that the problem (\ref{EDSM}) is well-posed. Their
proof is based on a fixed point theorem with respect to the
Kantorovitch-Rubeinstein metric.

\item At our knowledge, there are no existence and uniqueness results for (\ref{EDPM}) at least when $\Phi$ is not smooth.

\item By It\^o's formula, similarly to the proof of Proposition
\ref{PI1.3}, it is easy to see that a solution $Y^{\varepsilon}$ of
(\ref{EDSM}) provides a solution $v^{\varepsilon}$ of
 (\ref{EDPM}), in the sense of distributions.
\end{enumerate}
\end{remark}

When $\beta$ is non-degenerate it is possible to show that
formulations \eqref{EDPM} and \eqref{EDSM} are equivalent. In
particular we have the following result.
\begin{theorem} \label{TheoMollif}
We suppose that $\beta$ is non-degenerate and $\varepsilon>0$ is
fixed.
\begin{enumerate}
\item [(1)] If $Y^{\varepsilon}$ is a solution of \eqref{EDSM} then
$v^{\varepsilon}:[0,T]\rightarrow \mathcal{M}(\mathbb{R})$, where
$v^{\varepsilon}(t,\cdot)$ is the law of $Y^{\varepsilon}_t$, is a
solution of \eqref{EDPM} and fulfills the following property
\[\mbox{\textbf{(P)} $v^{\varepsilon}$ has
a density, still denoted $v^{\varepsilon}$ such that:}\ (t,x)\mapsto
v^{\varepsilon}(t,x) \in L^2([0,T]\times\mathbb{R}).\]

\item [(2)] If $v^{\varepsilon}$ is a solution to \eqref{EDPM}
fulfilling \textbf{(P)} then there is a process $Y=Y^{\varepsilon}$
solving \eqref{EDSM}.
\end{enumerate}
\end{theorem}

\begin{proof}

(1) If $Y^{\varepsilon}$ is a solution to \eqref{EDSM} by Remark
\ref{Remark_Mollif}.(iii) it follows that $v^{\varepsilon}$ fulfills
\eqref{EDPM}.

On the other hand, since $K_{\varepsilon}*v^{\varepsilon}$ is
bounded and $\Phi$ is lower bounded by a constant $C_{\varepsilon}$
on $[-\inf K_{\varepsilon}*v^{\varepsilon},\sup
K_{\varepsilon}*v^{\varepsilon}]$ it follows that
$a(t,x)=\Phi^2\left(K_{\varepsilon}*v^{\varepsilon}(t,x)\right)$ is
lower bounded by $C^{\varepsilon}$.

Using then Exercise 7.3.3 of \cite{StrVar}, i.e., Krylov estimates,
it follows that for every smooth function $f:[0,T]\times\mathbb{R}
\rightarrow \mathbb{R}$ with compact support, we have
\[E\left(\int_0^T f(Y_s^{\varepsilon})ds\right) \leq \mbox{const}\|f\|_{L^2([0,T]\times \mathbb{R})}.\]
Then, developing the left hand side to obtain
\[\int_0^T ds \int_{\mathbb{R}}f(y)v^{\varepsilon}(s,y)dy \leq \mbox{const}\|f\|_{L^2([0,T]\times \mathbb{R})}.\]
we deduce that (P) is verified.\\

(2) We retrieve here some arguments used in the proof of Proposition
4.2 of \cite{3}.

Given $v=v^{\varepsilon}$, by Remark 4.3
 of \cite{3}, see also
Exercise 7.3.2-7.3.4 of \cite{StrVar}, we can construct a  unique
solution $Y=Y^{\varepsilon}$ in law to the SDE constituted by
\begin{equation} \label{EDSmeas}
Y_t=Y_0+ \int\limits_0^t a(s,Y_s)dW_s,
\end{equation}
where here $a(t,x)=\Phi^2\left(K_{\varepsilon}*v(t,x)\right).$
Indeed, this is possible again because $a$ is Borel bounded and
lower bounded by a strictly postive constant.

A further use of It\^o's formula says that the law $z(t,dx)$ of
$Y_t$ solves
\begin{equation}\label{EDPZ}
 \left\{
 \begin{array}{ccl}
\partial_tz(t,.)&=&\frac{1}{2} \partial_{xx}^2 \left(a(t,.) z(t,.)\right),\\
       z(0,.)&=&u_0,\\
\end{array}
\right.
\end{equation}
in the sense of distributions.

Using again Krylov estimates as in the second part of the proof of
 point (1), it follows that $z$ admits a density $(t,y)\mapsto p_t(y)$ which
verifies $p \in L^2([0,T]\times \mathbb{R})$. This shows that
Hypothesis \ref{(A)} in Theorem \ref{T3.6} below is fulfilled, which
implies that $v\equiv z$.
\end{proof}
Theorem \ref{T3.6} was stated and proved in  \cite{3}, see Theorem
3.8.
\begin{theorem}\label{T3.6}
Let $a$ be a Borel nonnegative bounded function on $[0,T]\times
\mathbb{R}$.

Let $z_i:[0,T]\rightarrow \mathcal{M}_+(\mathbb{R})$, $i=1,2$, be
continuous with respect to the weak topology on finite measures on
$\mathcal{M}(\mathbb{R})$.

Let $z^0$ be an element of $\mathcal{M}_+(\mathbb{R})$. Suppose that
both $z_1$ and $z_2$ solve the problem $\partial_t
z=\partial_{xx}^2(az)$ in the sense of distributions with initial
condition $z(0,\cdot)=z^0$.

More precisely
\begin{equation*}
\int_{\mathbb{R}}\phi(x)z(t,dx)=\int_{\mathbb{R}}\phi(x)z^0(dx)+\int_0^t
ds \int_{\mathbb{R}}\phi''(x)a(s,x)z(s,dx)
\end{equation*}
for every $t\in[0,T]$ and any $\phi \in C_{0}^{\infty}(\mathbb{R})$.

Then $(z_1-z_2)(t,\cdot)$ is identically zero for every t, if
$z:=z_1-z_2$, satisfies the following:
\end{theorem}
\begin{hypothesis}\label{(A)}
  There is $\rho:[0,T]\times
\mathbb{R}\rightarrow \mathbb{R}$ belonging to $L^2([\kappa,T]\times
\mathbb{R})$ for every $\kappa >0$ such that $\rho(t,\cdot)$ is the
density of $z(t,\cdot)$ for almost all $t\in ]0,T]$.
\end{hypothesis}
\subsection{The interacting particles system}

We recall that in this paper, we want to approximate solutions of
 problem (\ref{EDP}). For this purpose we will concentrate on a probabilistic particles
system  of the same nature as in \cite{11} when the coefficients are
Lipschitz.

In  general, the particles probabilistic algorithms for non linear
PDEs are based on the  simulation of particles trajectories animated
by a random motion. The solution of the PDE
 is approximated
through the smoothing of the empirical measure of the particles,
which is
 a linear combination of Dirac masses centered on
particles positions. This procedure is heuristically justified by
the chaos propagation phenomenon which will be explained in the
sequel.

The dynamics of the particles is described by the following
stochastic differential system:
\begin{equation}\label{Sys(2)}
      Y_t^{i,\varepsilon,n}= Y_0^i+ \int_0^t \Phi\left(\frac{\displaystyle{1}}{\displaystyle{n}}\sum_{j=1}^n
K_{\varepsilon}(Y_s^{i,\varepsilon,n}-Y_s^{j,\varepsilon,n})\right)dW_s^i,\
i=1,\ldots,n
\end{equation}
where $W=(W^1,\ldots,W^n)$ is an n-dimensional Brownian motion,
$(Y_0^i)_{1\leq i\leq n}$ is a sequence of independent random
variables  with law density
 $u_0$ and independent of the Brownian motion $W$ and $K_\varepsilon$
is the same kernel as in Subsection \ref{Section: Some
Complements}.1.
\begin{remark}
If $\Phi$ in the system of  ordinary SDEs (\ref{Sys(2)}) were not
continuous but only measurable, that problem would not have
necessarily a solution, even if $\beta$ were non-degenerate. In
fact, contrarily to \eqref{EDSmeas}, here  $n \ge 2$. Since $\Phi$
is  continuous,  then \eqref{Sys(2)} has at least a solution; if
$\Phi$ is non-degenerate even uniqueness holds, see
 Chapter 6 and 7 of \cite{StrVar}.
\end{remark}
Now, owing to the interacting kernel $K_{\varepsilon }$, the
particles motions  are a priori dependent. For a given integer $n$,
we consider
 $(Y_t^{1,\varepsilon,n}, \ldots, Y_t^{n,\varepsilon,n})$ as the solution
of the interacting particle system  \eqref{Sys(2)}.
  Propagation of chaos  for the mollified equation happens if
for any integer $m$, the vector
 $(Y_t^{1,\varepsilon,n}, \ldots, Y_t^{m,\varepsilon,n})_{n \ge m}$
converges in law to $\mu_t \otimes^m$ where $\mu_t$ is the law of
$Y_t^{\varepsilon}$ the solution of \eqref{EDSM}.

A  consequence of chaos propagation is that one expects that
 the empirical measure of
the particles, i.e. the linear combination of Dirac masses denoted
$\mu^n_t=\frac{\displaystyle{1}}{\displaystyle{n}}\sum\limits_{j=1}^n
\delta_{Y_t^{j,\varepsilon,n}}$ converges in law as a random measure
to the deterministic solution $v_\varepsilon(t,.)$ of the
regularized PDE (\ref{EDPM}) which in fact depends on $\varepsilon$.
%Letting $\varepsilon$ go to zero $v_\varepsilon$.
This fact was established for instance when $\beta$ is Lipschitz, in
Proposition 2.2 of \cite{11}. On the other hand when $\varepsilon$
goes to zero, the same authors show that $v_\varepsilon$ converge to
the solution $u$ of \eqref{EDP}. They prove the existence of a
sequence $(\varepsilon(n))$ slowly converging to zero when $n$ goes
to infinity such  that the empirical measure
$\frac{\displaystyle{1}}{\displaystyle{n}}\sum\limits_{j=1}^n
\delta_{{Y_t^{j,\varepsilon(n),n}}}$, converges in law to $u$, see
Theorem 2.7 of \cite{11}. One consequence of the slow convergence is
that the regularized empirical measure
$$ \frac{\displaystyle{1}}{\displaystyle{n}}\sum\limits_{j=1}^n
K_{\varepsilon(n)}(\cdot-Y_t^{j,\varepsilon(n),n})$$ also converges
to $u$. Consequently, this probabilistic interpretation provides an
algorithm allowing to solve numerically  (\ref{EDP}).

  We recall however that one of the significant object of this paper is the numerical implementation related to the case when, $\beta$ is possibly discontinuous; for the moment we
 do not have convergence results but we implement the same
type of algorithm and we compare with some existing deterministic
schemes.
\section{About probabilistic numerical implementations}\label{Section:Num.
Impl.}

In this section we will try to construct an approximation method for
solutions $u$ of  \eqref{EDP}, based upon the time discretization of
the system (\ref{Sys(2)}). For now on, the number $n$ of particles
is fixed.

In fact, to get a simulation procedure for a trajectory of each
$(Y_t^{i,\varepsilon,n}),~~i=1,\ldots,n$, we discretize in time: for
fixed $T>0$, we choose  a time step $\Delta t>0$ and $N \in
\mathbb{N}$, such that $T=N\Delta t$. We denote by $t_k=k \Delta t$,
the discretization times for $k=0,\ldots,N$.

The Euler  explicit scheme of order one, leads then to the following
discrete time system, i.e., for every $ i=1,\ldots,n$
\begin{eqnarray}\label{SysDiscr}
      \begin{aligned}
X_{t_{k+1}}^i &=& X_{t_k}^i +\Phi
\left(\frac{\displaystyle{1}}{\displaystyle{n}}\sum_{j=1}^n
K_{\varepsilon} ( X_{t_k}^i- X_{t_k}^j)
\right)\left(W_{t_{k+1}}^i-W_{t_k}^i\right),
      \end{aligned}
\end{eqnarray}
where at each time step $t_k$, we approximate $u(t_k,.)$ by the
smoothed empirical measure of the particles :
\begin{equation}\label{estimateur}
u^{\varepsilon,n}(t_k,x)=\frac{\displaystyle{1}}{\displaystyle{n}}\sum_{j=1}^n
K_{\varepsilon}(x-X_{t_k}^j),~~k=1,\ldots,N,\ \ \ \ x \in
\mathbb{R},
\end{equation}
at each time step and for every $i=1,\ldots,n$, the Brownian
increment $\left(W_{t_{k+1}}^i-W_{t_k}^i\right)$ is given by the
simulation of the realization of a Gaussian random variable of law
$\mathcal{N}(0,\Delta t)$.

One difficult issue concerns the smoothing parameter $\varepsilon$
related to the kernel $K_{\varepsilon}$. It will be chosen according
to the \emph{kernel density estimation}.

In fact from now on we will assume that $K$, as defined in
(\ref{Scaling}), is a Gaussian probability density function with
mean $0$ and unit standard deviation. In this case,
 in (\ref{estimateur}), the function
$u^{\varepsilon,n}(t_k,\cdot)$  becomes the so-called Gaussian
kernel density estimator of $u(t_k,\cdot)$  for every time step
$t_k$ with $k=1,\ldots,N$.

Finally, the only unknown parameter in (\ref{estimateur}), is
 $\varepsilon$; most of the authors  refer to it as the \emph{bandwidth}
or the \emph{window width}.

The optimal choice of $\varepsilon$ was the object of an enormous
amount of research, because its value strongly determines the
performance of $u^{\varepsilon,n}$ as an estimator of $u$ depends,
see, e.g. \cite{14} and references therein. The most widely used
criterion of performance for the estimator (\ref{estimateur}) is the
\emph{Mean Integrated Squared Error} (MISE), defined by
\begin{eqnarray*}
\mbox{MISE}\{u^{\varepsilon,n}(t,x)\}&=&\mathbb{E}_{u}\int\left[u^{\varepsilon,n}(t,y)-u(t,y)\right]^2
dy \\
&=&\int
\left(\underbrace{\mathbb{E}_{u}\left[u^{\varepsilon,n}(t,y)\right]-u(t,y)}_{\mbox{\footnotesize{point-wise
bias}}}\right)^2dy +\underbrace{\int
\mathbb{V}_u\left[u^{\varepsilon,n}(t,y)\right]dy,}_{\mbox{\footnotesize{
integrated point-wise variance}}}
\end{eqnarray*}
where, $\mathbb{E}_u$ and $\mathbb{V}_u$ are respectively the
expectation and the variance of $X^j_t,\ j=1,..,n$ under the
assumption that they are independent and distributed as
$u(t,\cdot)$.

We emphasize that the MISE expression is the sum of two components:
the integrated bias and  variance.

The asymptotic properties of (\ref{estimateur}) under the MISE
criterion are well-known (see\cite{14},\cite{BotevPreparation}), but
we summarize them below for convenience of the reader.

\begin{theorem}\textbf{(Properties of the Gaussian kernel estimator)}

Under the assumption that $\varepsilon$ depends on $n$ such that
$\lim\limits_{n\to {+\infty}}\varepsilon=0,$ $\lim\limits_{n \to
{+\infty}} n\varepsilon=+\infty $ and $\partial^2_{xx}u$ is a
continuous square integrable function, the estimator
(\ref{estimateur}) has integrated squared bias and integrated
variance given by
\begin{eqnarray}
\|E_u\left[u^{\varepsilon,n}(t,.)-u(t,.)\right]\|^2&=&\frac{\displaystyle{1}}{\displaystyle{4}}\varepsilon^4
\|\partial^2_{xx}u\|^2 + o(\varepsilon^2),~~n\rightarrow +\infty,\\
\int
V_u\left[u^{\varepsilon,n}(t,y)\right]dy&=&\frac{\displaystyle{1}}{\displaystyle{2\varepsilon
n \sqrt{\pi}}}+o((n \varepsilon)^{-1}),~~n\rightarrow +\infty.
\end{eqnarray}
\end{theorem}

\begin{remark}
\begin{enumerate}
    \item  Here $\|.\|$ denotes the standard $L^2$ norm. The first order
asymptotic approximation of MISE, denoted AMISE, is thus given by
\begin{equation}\label{AMISE}
\mbox{AMISE}\{u^{\varepsilon,n}(t,x)\}=\frac{\displaystyle{1}}{\displaystyle{4}}\varepsilon^4
\|\partial^2_{xx}u(t,x)\|^2 + (2\varepsilon n \sqrt{\pi})^{-1}.
\end{equation}

    \item  The asymptotically optimal value of $\varepsilon$ is the minimizer
of AMISE and by simple calculus it can be shown (see \cite{15},
Lemma 4A) to be equal to $\varepsilon_t^{opt}$ defined in formula
\eqref{OptimalEpsilon}.
\end{enumerate}

\end{remark}

As argued in the introduction, we have chosen to use the
"solve-the-equation" bandwidth selection plug-in procedure developed
in \cite{SJ,Jones_al_1996}, to perform the optimal window width of
the Gaussian kernel density estimator $u^{\varepsilon,n}$ of $u$,
defined in \eqref{estimateur}.

\begin{remark}
According to \cite{SJ}, for every positive integer $s$, the identity
\[\|\partial^s_{x^s}u(t,.)\|^2=(-1)^s\int_{\mathbb{R}}\partial^{2s}_{x^{2s}}u(t,x).u(t,x)dx,\]
suggests the following estimator for that density functional:
\begin{equation}\label{SJ_estimator}
    \|\partial^s_{x^s}u^{\varepsilon,n}(t,x)\|^2=\frac{\displaystyle{(-1)^s}}{\displaystyle{n^2\varepsilon^{2s+1}}}
  \sum_{i=1}^n \sum_{j=1}^n K^{(2s)}\left(\frac{\displaystyle{X_t^i-X_t^j}}{\displaystyle{\varepsilon}}\right)
\end{equation}
where, $K^{(r)}$ is the $r^{th}$ derivative of the Gaussian kernel
$K$ and $\partial^r_{x^r}u$ is the $r^{th}$ partial spacial
derivative of $u$.
\end{remark}

Inspired, by \eqref{OptimalEpsilon}, the authors of
\cite{SJ,Jones_al_1996} look for an approached  optimal bandwidth
for the AMISE as the solution of the equation
\begin{equation}\label{AMISE-OptEps}
\varepsilon:=\varepsilon_t={\left(2n\sqrt{\pi}
\|\partial^2_{xx}u^{\gamma(\varepsilon_t),n}(t,x)\|^2\right)}^{-1/5}
\end{equation}
where, $\|\partial^2_{xx}u^{\gamma(\varepsilon),n}\|^2$ is an
estimate of $\|\partial^2_{xx}u\|^2$ using \eqref{SJ_estimator}, for
$s=2$, and the pilot bandwidth $\gamma(\varepsilon)$, which depends
on the kernel bandwidth $\varepsilon$. The pilot bandwidth
$\gamma(\varepsilon)$ is  then chosen through an intermediate step
which consists in obtaining a quantity $h_t$
 minimizing the
asymptotic mean squared error (AMSE) for the estimation of
$\|\partial^2_{xx}u\|^2$. AMSE is in fact some approximation via
Taylor expansion of the MSE which is defined as follows:
\begin{equation}\label{AMSE}
    \mbox{MSE}\{\|\partial^2_{xx}u^{h,n}\|^2\}=\mathbb{E}_u{\left[\|\partial^2_{xx}u^{h,n}\|^2-\|\partial^2_{xx}u\|^2\right]}^2.
\end{equation}
Similarly one can  define an analogous quantity for the third
derivative
\begin{equation}\label{AMSE*}
    \mbox{MSE}\{\|\partial^3_{x^3}u^{h,n}\|^2\}=\mathbb{E}_u{\left[\|\partial^3_{x^3}u^{h,n}\|^2-\|\partial^3_{x^3}u\|^2\right]}^2.
\end{equation}
and related AMSE. Exhaustive details concerning those computations
are given  in \cite{Wand_Jones_1995}. In fact, the authors in
\cite{Wand_Jones_1995} computed those minimizers and provided the
following explicit formulae
\begin{equation}\label{h_pilotEps}
    h_t={\left[\frac{\displaystyle{2K^{(4)}(0)}}{\displaystyle{n
    \|\partial^3_{x^3}u(t,x)\|^2}}\right]}^{1/7},   \quad
 h^*_t={\left[\frac{\displaystyle{-2K^{(6)}(0)}}{\displaystyle{n
    \|\partial^4_{x^4}u(t,x)\|^2}}\right]}^{1/9},
\end{equation}
where $h_t$ and    $h^*_t$ minimize the AMSE corresponding
respectively to \eqref{AMSE} and \eqref{AMSE*}.

Solving \eqref{OptimalEpsilon}, with respect to $n$ and replacing
$n$ in the first equality of \eqref{h_pilotEps}, gives the following
expression of $h_t$ in term of $\varepsilon_t$
\begin{equation*}
    h_t={\left[\frac{\displaystyle{4\sqrt{\pi}K^{(4)}(0)\|
\partial^2_{xx}u(t,x)\|^2}}
    {\displaystyle{
    \|\partial^3_{x^3}u(t,x)\|^2}}\right]}^{1/7}\varepsilon^{5/7}_t.
\end{equation*}
This suggests to define
\begin{equation}\label{Gamma_h}
\gamma(\varepsilon_t)={\left[\frac{\displaystyle{4\sqrt{\pi}
K^{(4)}(0)\|\partial^2_{xx}u^{h^1_t,n}(t,x)\|^2}}
    {\displaystyle{
    \|\partial^3_{x^3}u^{h^2_t,n}(t,x)\|^2}}\right]}^{1/7}\varepsilon^{5/7}_t,
\end{equation}
where, $\|\partial^2_{xx}u^{h^1_t,n}\|^2$ and
$\|\partial^3_{x^3}u^{h^2_t,n}\|^2$ are estimators of
$\|\partial^2_{xx}u\|^2 $ and $\|\partial^3_{x^3}u\|^2$ using
formula \eqref{SJ_estimator} and
 pilot
bandwidths $h_t^1$ and $h_t^2$ given by
\[h_t^1={\left[\frac{\displaystyle{2K^{(4)}(0)}}{\displaystyle{n \widehat{\|\partial^3_{x^3}u(t,x)\|^2}}}\right]}^{1/7} \ \ \ \ \
h_t^2={\left[\frac{\displaystyle{-2K^{(6)}(0)}}{\displaystyle{n
\widehat{\|\partial^4_{x^4}u(t,x)\|^2}}}\right]}^{1/9}\ ;\]
$\widehat{\|\partial^3_{x^3}u(t,x)\|^2}$ and
$\widehat{\|\partial^4_{x^4}u(t,x)\|^2}$ will be suitably defined
below. Indeed, $h_t^1$ and $h_t^2$ estimate $h_t$ and  $h^*_t$
defined in \eqref{h_pilotEps}.

According to the strategy in \cite{SJ, Wand_Jones_1995}, we will
first suppose  that $\partial^3_{x^3}u(t,x)$ and
$\partial^4_{x^4}u(t,x)$ are the third and fourth partial space
derivatives of a Gaussian density with standard deviation $\sigma_t$
of  $X_t$. In a second step we replace $\sigma_t$ with the empirical
 standard deviation $\hat{\sigma}_t$
 of the sample $X_t^1,\ldots,X_t^n$. This leads naturally  to
\[\widehat{\|\partial^3_{x^3}u(t,x)\|^2}=\frac{\displaystyle{15}}{\displaystyle{16
\sqrt{\pi}}}\hat{\sigma}^{-7}_t, \ \ \ \ \ \ \
\widehat{\|\partial^4_{x^4}u(t,x)\|^2}=\frac{\displaystyle{105}}{\displaystyle{32
\sqrt{\pi}}}\hat{\sigma}^{-9}_t.\]
 Coming back to
  \eqref{AMISE-OptEps}, where $\gamma(\varepsilon_t)$ is defined through
\eqref{Gamma_h}, it suffices then to perform a root-finding
algorithm for it  at each discrete time step $t_k$, in order to
obtain the approached optimal bandwidth $\varepsilon_{t_k}$.
\section{Deterministic numerical approach}\label{Section:Det. App}

We recall that the final aim of our work is to approximate solutions
of  a nonlinear problem given by
\begin{equation}\label{ProbD}
 \left\{
 \begin{array}{ccl}
  \partial_tu(t,x)&\in &\frac{1}{2} \partial_{xx}^2 \beta\left(u(t,x)\right),\ \ t\in\left[0,+\infty\right[,\\
        u(0,x)&=&u_0(x), \ \  \ x \in \mathbb{R},\\
\end{array}
\right.
\end{equation}
in the case where $\beta$ is given by (\ref{Heav}). Despite the fact
that, up to now at our knowledge,
 there are no analytical approaches
dealing such issues, we  got interested into a recent method,
proposed in \cite{17}. Actually, we are heavily inspired by
\cite{17} to implement a deterministic procedure  simulating
solutions of \eqref{ProbD} which will be compared to the
probabilistic one. \cite{17} handles with the propagation of a
discontinuous solution, even though coefficient $\beta$ is
Lipschitz. It seems to us that in the numerical analysis literature,
\cite{17} is the closest one to our spirit. We describe now the
fully discrete scheme we will use for this purpose.
\subsection{Relaxation approximation}
The schemes proposed  in \cite{17} follow the same idea as the
well-known relaxation schemes for hyperbolic conservation laws, see
\cite{Jin} for a review of the subject. For the convenience of the
reader, we retrieve here some arguments of \cite{17}, where we
recall that the coefficient $\beta$ is Lipschitz. In that case
$\in$, of course, becomes $=$.

The equation \eqref{ProbD} can be formally expressed by the first
order system on $\mathbb{R}_+\times \mathbb{R}$ :
\begin{equation} \label{RS1_Brut}
 \begin{cases}
\partial_t u+ \partial_x v=0, \\
v+ \frac{\displaystyle{1}}{\displaystyle{2}}\partial_x\beta(u)=0.
\end{cases}
\end{equation}
\eqref{RS1_Brut}, is relaxed with the help of a  parameter
$\varepsilon>0$, in order to obtain the following  scheme
\begin{equation} \label{RS1}
 \begin{cases}
\partial_t u+ \partial_x v=0, \\
\partial_t v+ \frac{\displaystyle{1}}{\displaystyle{2\varepsilon}}\partial_x\beta(u)=
-\frac{\displaystyle{1}}{\displaystyle{\varepsilon}}v.
\end{cases}
\end{equation}
Then, another function $w:\mathbb{R}_+ \times \mathbb{R} \rightarrow
\mathbb{R}$ is introduced in order to remove the non-linear term in
the second line of system \eqref{RS1}. So, we obtain
\begin{equation} \label{RS1Bis}
 \begin{cases}
\partial_t u+ \partial_x v=0, \\
\partial_t v+ \frac{\displaystyle{1}}{\displaystyle{2\varepsilon}}\partial_xw=
-\frac{\displaystyle{1}}{\displaystyle{\varepsilon}}v,\\
\partial_t w+ \partial_x
v=-\frac{\displaystyle{1}}{\displaystyle{\varepsilon}}(w-\beta(u)).
\end{cases}
\end{equation}
Note that \eqref{RS1Bis} is a particular case of the BGK system
previously studied in \cite{BKGSyst}. In fact, authors of
\cite{BKGSyst} proved that $w$ (resp. $v$) converges to $\beta(u)$
(resp.  $-\frac{\displaystyle{1}}{\displaystyle{2}}
\partial_x \beta(u)$), as  $\varepsilon \rightarrow 0^+$.
Furthermore, they showed the convergence of solutions of
\eqref{RS1Bis} to those of PDE \eqref{ProbD}, in $L^1(\mathbb{R})$,
as $\varepsilon$ goes to zero.

 Finally, we introduce a supplementary parameter $\varphi>0$,
 according to usual numerical analysis techniques; while preserving the hyperbolic character of the system.
 Therefore, we get
\begin{equation} \label{RS}
 \begin{cases}
\partial_t u+ \partial_x v=0, \\
\partial_t v+
\varphi^2\partial_x w=
-\frac{\displaystyle{1}}{\displaystyle{\varepsilon}}v +
(\varphi^2-\frac{\displaystyle{1}}{\displaystyle{2\varepsilon}})\partial_x
w ,\\
\partial_t w+ \partial_x v=
-\frac{\displaystyle{1}}{\displaystyle{\varepsilon}}(w-\beta(u)).
\end{cases}
\end{equation}
 Now, setting
\[z=\left (
\begin{array}{ccc}
   u \\
   v \\
   w
\end{array}
\right),\ \ \mathcal{F}(z)=\mathbb{A}z,\ \ \mathbb{A}=\left(
              \begin{array}{ccc}
                0 & 1 & 0 \\
                0 & 0 & \varphi^2 \\
                0 & 1 & 0 \\
              \end{array}
            \right) \ \ \mbox{and}\ \ g(z)=\left (
\begin{array}{ccc}
   0 \\
   -v + (\varphi^2 \varepsilon -\frac{1}{2})\partial_x w \\
   \beta(u)-w
\end{array}
\right),\] the system (\ref{RS}) is rewritten in  matrix form as
follows
\begin{equation}\label{RSM}
\partial_t z+
\partial_x \mathcal{F}(z)= \frac{\displaystyle{1}}{\displaystyle{\varepsilon}}
g(z).
\end{equation}
Using the change of variable $Z=\mathbb{P}^{-1}z$, where
\[\mathbb{P}^{-1}=\left(
                         \begin{array}{ccc}
                           0 & \frac{\displaystyle{1}}{\displaystyle{2 \varphi}} & \frac{\displaystyle{1}}{\displaystyle{2}} \\
                           0 & \frac{\displaystyle{-1}}{\displaystyle{2 \varphi}} & \frac{\displaystyle{1}}{\displaystyle{2}}\\
                           1 & 0 & -1 \\
                         \end{array}
                       \right)\  \mbox{ and } \ \mathbb{P}^{-1} \mathbb{A} \mathbb{P}=\mathbb{D}=\left(
                                                          \begin{array}{ccc}
                                                            \varphi & 0 & 0 \\
                                                            0 & -\varphi & 0\\
                                                            0 & 0 & 0 \\
                                                          \end{array}
                                                        \right),\]
 we obtain
\begin{equation}\label{Rel1}
Z=\left (
\begin{array}{ccc}
   \mathcal{U}\\
   \mathcal{V}\\
\mathcal{W}
\end{array}
\right),\ \ \mbox{with}\ \ \ \mathcal{U}=\frac{\displaystyle{ v
+\varphi w}}{\displaystyle{2\varphi}},\ \ \mathcal{V}
=\frac{\displaystyle{ -v +\varphi w}}{\displaystyle{2\varphi}},\ \
\mathcal{W}=u-w,
\end{equation}
 where,  $ \ \mathcal{U}$,  $\mathcal{V}$, $\mathcal{W}$ are called
characteristic variables.

Since $z=\mathbb{P} Z$, equation (\ref{RSM}) leads to
\begin{equation}\label{RSMI}
\partial_t Z+
\mathbb{D}  \partial_x
Z=\frac{\displaystyle{1}}{\displaystyle{\varepsilon}}
\mathbb{P}^{-1}g(\mathbb{P}Z).
\end{equation}
By rewriting the system (\ref{RSMI}) in terms of the characteristic
variables, we obtain
\begin{equation}\label{SUVW}
\left (
\begin{array}{ccc}
   \partial_t
\mathcal{U}+ \varphi \partial_x \mathcal{U}\\ \\
   \partial_t{V}-\varphi
\partial_x \mathcal{V}\\ \\
\partial_t \mathcal{W}
\end{array}
\right)=\frac{\displaystyle{1}}{\displaystyle{\varepsilon}}
\mathbb{P}^{-1}g(\mathbb{P}Z).
\end{equation}
Finally,  solving \eqref{RS} is equivalent to the resolution of a
three advection equations system, \eqref{SUVW}, with respectively a
positive, a negative and a zero advection velocity.

\begin{remark}\label{Remark_uUVW}
Note that we can deduce from \eqref{Rel1} the following relation
\begin{equation}\label{Rel2}
u=\mathcal{U}+\mathcal{V}+\mathcal{W}.
\end{equation}
\end{remark}
\subsection{Space discretization}
In the sequel of this chapter and in Annex \ref{Annex}, given two
integers $i < j $,
$\left[\hspace{-1ex}\left[\hspace{0.5ex}i,j\hspace{0.5ex}
\right]\hspace{-1ex}\right ]$ , will denote the integer interval
$\{i, i+1,\ldots, j \}$. We will now provide a space discretization
scheme for  system \eqref{SUVW}. Let us introduce a uniform grid on
$[a,b] \subset \mathbb{R}$.

 We denote
$x_i=a-\frac{\displaystyle{\Delta x}}{\displaystyle{2}}+i\Delta x~$,
$ i \in \left[\hspace{-1ex}\left[\hspace{0.5ex}1,N_x\hspace{0.5ex}
\right]\hspace{-1ex}\right]\ $  and  $\ x_{i+1/2}=a+i\Delta x~$, $ i
\in \left[\hspace{-1ex}\left[\hspace{0.5ex}0,N_x\hspace{0.5ex}
\right]\hspace{-1ex}\right] \ $, where $\Delta x=
\frac{\displaystyle{b-a}}{\displaystyle{N_x}}$ is the grid spacing
and $N_x$ the number of cells. Note that $x_i$ is the center of the
interval $[x_{i-1/2},x_{i+1/2}]$. Moreover, we denote the boundary
conditions by $u(t,a)=u_a(t)$ and $u(t,b)=u_b(t)$, for every $t>0$.

Then, we evaluate \eqref{SUVW} on the grid of discrete points
$(x_i)$ getting,
\begin{equation}\label{DS2}
 \left\{
 \begin{array}{rcl}
  \frac{\displaystyle{d\mathcal{U}}}{\displaystyle{dt}}(t,x_i)+
\varphi
\frac{\displaystyle{d\mathcal{U}}}{\displaystyle{dx}}(t,x_i)&=&G_1(t,x_i),\
\ \forall t>0,\
\forall i \in \left[\hspace{-1ex}\left[\hspace{0.5ex}1,N_x\hspace{0.5ex} \right]\hspace{-1ex}\right],\\ \\
\frac{\displaystyle{d\mathcal{V}}}{\displaystyle{dt}}(t,x_i)-
\varphi
\frac{\displaystyle{d\mathcal{V}}}{\displaystyle{dx}}(t,x_i)&=&G_2(t,x_i),\
\  \forall t>0, \
\forall i \in \left[\hspace{-1ex}\left[\hspace{0.5ex}1,N_x\hspace{0.5ex} \right]\hspace{-1ex}\right],\\ \\
\frac{\displaystyle{d\mathcal{W}}}{\displaystyle{dt}}(t,x_i)&=&G_3(t,x_i),\
\ \forall t>0,\ \forall i \in
\left[\hspace{-1ex}\left[\hspace{0.5ex}1,N_x\hspace{0.5ex}
\right]\hspace{-1ex}\right],
\end{array}
\right.
\end{equation}
where,
\begin{equation}\label{G1_G2_G3}
(G_1,G_2,G_3)^t=\frac{\displaystyle{1}}{\displaystyle{\varepsilon}}
\mathbb{P}^{-1}g(\mathbb{P}Z).
\end{equation}
\begin{remark}\label{Somme_Gi}
We can easily deduce from \eqref{G1_G2_G3},  that for every $t \in
]0,+\infty[\ $ and every $\ i \in
\left[\hspace{-1ex}\left[\hspace{0.5ex}1,N_x\hspace{0.5ex}
\right]\hspace{-1ex}\right]$ , we have : $\ \sum\limits_{j=1}^3
G_j(t,x_i)=0.$
\end{remark}
In order to ensure the convergence of the semi-discrete scheme
\eqref{DS2} it is necessary to write it in a conservative form. To
this aim, following \cite{17},  we suppose the existence of
functions $\widehat{\mathcal{U}}$ and $\widehat{\mathcal{V}}$ such
that
\begin{equation*}
\begin{cases}
\mathcal{U}(t,x)=\frac{\displaystyle{1}}{\displaystyle{\Delta x}
}\int\limits_{x-\Delta x/2}^{x+\Delta x/2}\widehat{\mathcal{U}}(t,y)
dy,\ \ \forall x \in ]a,b[,\ \ \forall t
>0,\\
\mathcal{V}(t,x)=\frac{\displaystyle{1}}{\displaystyle{\Delta x} }
\int\limits_{x-\Delta x/2}^{x+\Delta x/2}\widehat{\mathcal{V}}(t,y)
dy,\ \ \forall x \in ]a,b[,\ \ \forall t
>0.
\end{cases}
\end{equation*}
Substituting in \eqref{DS2}, we obtain  for every $ t>0\ $ and every
$\ i \in \left[\hspace{-1ex}\left[\hspace{0.5ex}1,N_x\hspace{0.5ex}
\right]\hspace{-1ex}\right]$ ,
\begin{equation}\label{DS2-3}
\left\{
 \begin{array}{rcl}
 \frac{\displaystyle{d\mathcal{U}}}{\displaystyle{dt}}(t,x_i)+
 \frac{\displaystyle{\varphi}}{\displaystyle{\Delta x}}
\left(\widehat{\mathcal{U}}(t,x_{i+1/2})-\widehat{\mathcal{U}}(t,x_{i-1/2})
\right)&=&G_1(t,x_i),\\ \\
\frac{\displaystyle{d\mathcal{V}}}{\displaystyle{dt}}(t,x_i)-
\frac{\displaystyle{\varphi }}{\displaystyle{\Delta x}}
\left(\widehat{\mathcal{V}}(t,x_{i+1/2})-\widehat{\mathcal{V}}(t,x_{i-1/2})
\right)&=&G_2(t,x_i),\\ \\
\frac{\displaystyle{d\mathcal{W}}}{\displaystyle{dt}}(t,x_i)&=&G_3(t,x_i).
\end{array}
\right.
\end{equation}
Let us now denote by $\ \widetilde{\mathcal{U}}_{i+1/2}(t)\ $ and $\
\widetilde{\mathcal{V}}_{i+1/2}(t)\ $ the so-called {\it
semi-discrete numerical fluxes} that approximate respectively $\
\widehat{\mathcal{U}}(t,x_{i+1/2})$ and $\
\widehat{\mathcal{V}}(t,x_{i+1/2})$. For the sake of simplicity, we
chose to expose only the calculations necessary to obtain the first
semi-discrete flux $\ \widetilde{\mathcal{U}}_{i+1/2}(t)$, the same
procedure being applied for the other one.

In order to compute the numerical flux $\
\widetilde{\mathcal{U}}_{i+1/2}(t)$, we reconstruct boundary
extrapolated data $~\mathcal{U}_{i+1/2}^{\pm}(t)$, from the point
values $~\mathcal{U}_i(t)=\mathcal{U}(t,x_i)$ of the variables at
the center of the cells, with an essentially non oscillatory
interpolation (ENO) method. The ENO technique allows to better
localize discontinuities and fronts that may appear when $\beta$ is
possibly degenerate; see \cite{18,Shu} for an extensive presentation
of the subject. In fact, $\ \mathcal{U}_{i+1/2}^{+}(t)$  (resp. $\
\mathcal{U}_{i+1/2}^{-}(t)$) is calculated from an interpolating
polynomial of degree $d$, on the interval $[x_{i+1/2},x_{i+3/2}]$
(resp. $[x_{i-1/2},x_{i+1/2}]$) using a so-called \emph{ENO
stencil}, see \cite{Shu} and formula \eqref{InterpData2} in Annex
\ref{SecAnnex1}.

Next, we shall apply a numerical flux to these boundary extrapolated
data. In order to minimize the numerical viscosity and according to
authors of \cite{17}, we choose the so-called \emph{Godunov flux},
$\mathfrak{F}_G $, associated to the advection equation
\begin{equation*}
    \partial_t \mathcal{U}+ \partial_x f(\mathcal{U})=0,
\end{equation*}
and defined as follows
\begin{equation*}
\mathfrak{F}_G [\alpha,\gamma] =\left\{\begin{array}{ll}
\min\limits_{\alpha\leq \xi \leq \gamma} f(\xi),\ \ & \mbox{if }\ \alpha\leq \gamma,\\ \\
\max\limits_{\gamma\leq \xi \leq \alpha} f(\xi),\ \ & \mbox{if }\ \gamma\leq \alpha.\\
\end{array}\right.
\end{equation*}
where $f(\xi)=\varphi \xi$, with $\varphi>0$. So we have,
$\mathfrak{F}_G [\alpha,\gamma] =\varphi\alpha$.

 In fact, we set
\begin{equation}\label{def1}
\forall t>0, \ \ \ \ \ \
\widetilde{\mathcal{U}}_{i+1/2}(t)=\mathfrak{F}_G[\mathcal{U}_{i+1/2}^{-}(t),\mathcal{U}_{i+1/2}^{+}(t)].
\end{equation}
Therefore, we obtain the following semi-discrete flux
\begin{equation}\label{Fluxes}
\forall t>0, \ \ \ \
\widetilde{\mathcal{U}}_{i+1/2}(t)=\varphi\mathcal{U}_{i+1/2}^{-}(t).
\end{equation}
 Applying the previous procedure to compute $\widetilde{\mathcal{V}}_{i+1/2}(t)$ and replacing  in \eqref{DS2-3}, we
 get for every $t>0\ $ and every $\ i \in
\left[\hspace{-1ex}\left[\hspace{0.5ex}1,N_x\hspace{0.5ex}
\right]\hspace{-1ex}\right]$ ,
\begin{equation}\label{DS3}
\left\{
 \begin{array}{rcl}
 \frac{\displaystyle{d\mathcal{U}}}{\displaystyle{dt}}(t,x_i)+
\frac{\displaystyle{\varphi}}{\displaystyle{\Delta x}}
\left(\mathcal{{U}}_{i+1/2}^{-}(t)-\mathcal{{U}}_{i-1/2}^{-}(t)
\right)&=& G_1(t,x_i),\\ \\
\frac{\displaystyle{d\mathcal{V}}}{\displaystyle{dt}}(t,x_i)-
\frac{\displaystyle{\varphi}}{\displaystyle{\Delta x}}
\left(\mathcal{{V}}_{i+1/2}^{+}(t)-\mathcal{{V}}_{i-1/2}^{+}(t)
\right)&=& G_2(t,x_i), \\ \\
\frac{\displaystyle{d\mathcal{W}}}{\displaystyle{dt}}(t,x_i)&=&G_3(t,x_i).
\end{array}
\right.
\end{equation}
Consequently summing up the three equation lines in \eqref{DS3} and
using Remarks \ref{Remark_uUVW} and \ref{Somme_Gi}, we obtain
\begin{equation*}
\frac{\displaystyle{du}}{\displaystyle{dt}}(t,x_i)+
\frac{\displaystyle{\varphi}}{\displaystyle{\Delta
x}}\left(\mathcal{{U}}_{i+1/2}^{-}(t)-
\mathcal{{U}}_{i-1/2}^{-}(t)-\left(\mathcal{{V}}_{i+1/2}^{+}(t)-\mathcal{{V}}_{i-1/2}^{+}(t)\right)
\right)=0.
\end{equation*}
Now, coming back to the conservative variables, we obtain for every
$i\in \left[\hspace{-1ex}\left[\hspace{0.5ex}1,N_x\hspace{0.5ex}
\right]\hspace{-1ex}\right]$ and every $t>0$,
\begin{eqnarray}\label{SysEDO}
\begin{cases}
 \frac{\displaystyle{du}}{\displaystyle{dt}}(t,x_i)=
-\frac{\displaystyle{1}}{\displaystyle{2\Delta
x}}\left(v_{i+1/2}^-(t)- v_{i-1/2}^-(t)+\varphi (w_{i+1/2}^-(t)-
w_{i-1/2}^-(t))\right)\\
\hspace{11.5ex}+\frac{\displaystyle{1}}{\displaystyle{2\Delta
x}}\left(v_{i-1/2}^+(t)-v_{i+1/2}^+(t)+\varphi
(w^+_{i+1/2}(t)-w^+_{i-1/2}(t)) \right),\\
u(0,x_i)\ \  = \ u_0(x_i),\\
u(t,a)\ \  \ \  = \ u_a(t),\\
u(t,b)\ \  \ \  = \ u_b(t).
\end{cases}
\end{eqnarray}
 We recall that by formally setting $\varepsilon=0$ in the scheme \eqref{RS}, we have $v=-\frac{\displaystyle{1}}{\displaystyle{2}}\partial_x w$ and $w=\beta(u)$.
 Therefore we can compute
\[v_{i+1/2}^{\pm}=-\frac{1}{2}
(\partial_x w)_{i+1/2}^{\pm} \ \ \mbox{and} \ \
w_{i+1/2}^{\pm}=\beta(u_{i+1/2}^{\pm}),\ \
\]
where, $w_{i+1/2}^{\pm}$, are performed using again an ENO
reconstruction, see formulae \eqref{InterpData1}-\eqref{InterpData2}
in Annex \ref{SecAnnex1}; while the derivatives of $w_{i+1/2}^{\pm}$
are approximated using a reconstruction polynomial with a centered
stencil, see formula
\eqref{DerivInterpData1}-\eqref{DerivInterpData3} in Annex
\ref{SecAnnex2}.

We wish to emphasize that the scheme of system \eqref{SUVW} reduces
to the time advancement of the single variable $u$ solution of
\eqref{ProbD}.
\subsection{Time discretization}

In order to have a fully discrete scheme, we  still need to specify
the time discretization. According to \cite{17}, we use a
discretization based on an  explicit Runge-Kutta  scheme, see
\cite{ParRusso},  for instance.

We start discretizing the system (\ref{SysEDO}) using,  for
simplicity, a uniform time step $\Delta t$. For every $i \in
\left[\hspace{-1ex}\left[\hspace{0.5ex}1,N_x\hspace{0.5ex} \right]
\hspace{-1ex}\right]$, we denote by $u_i^m$ the numerical
approximation of  $u(t^m,x_i)$ with $t^m=m\Delta t,\ \
m=0,\ldots,N_t$, where $N_t$ is the number of time steps.

The $\nu$-stage explicit Runge-Kutta scheme with $\nu \geq 1$,
associated to (\ref{SysEDO}) can be written  for every $i \in
\left[\hspace{-1ex}\left[\hspace{0.5ex}1,N_x\hspace{0.5ex} \right]
\hspace{-1ex}\right]$,  as follows,
\begin{equation}
u_i^{m+1}=u_i^m - \frac{\displaystyle{\lambda}}{\displaystyle{2}}
\sum_{k=1}^{\nu} \tilde{b}_k F_i^{(k)},
\end{equation}
where, $\lambda=\frac{\displaystyle{\Delta t}}{\displaystyle{\Delta
x}}$ and the stage values are computed at each time step $t^m$ and
for every $k \in
\left[\hspace{-1ex}\left[\hspace{0.5ex}1,\nu\hspace{0.5ex}
\right]\hspace{-1ex}\right]$, as
\begin{eqnarray}
\begin{cases}
F^{(k)}_i= v_{i+1/2}^{(k)-}- v_{i-1/2}^{(k)-}+\varphi
(w_{i+1/2}^{(k)-}-
w_{i-1/2}^{(k)-})-v_{i-1/2}^{(k)+}+v_{i+1/2}^{(k)+}-\varphi
(w^{(k)+}_{i+1/2}-w^{(k)+}_{i-1/2}),\\ \\
u_i^{(k)}=u_i^m-\frac{\displaystyle{\lambda}}{\displaystyle{2}}
\sum\limits_{l=1}^{k-1}\tilde{a}_{kl}F_i^{(l)},\ \
v_{i+1/2}^{(l)\pm}=-\frac{1}{2} (\partial_x w^{(l)})_{i+1/2}^{\pm},
\ \ w_{i+1/2}^{(l)\pm}=\beta(u_{i+1/2}^{(l)\pm}). \label{stage_uk}
\end{cases}
\end{eqnarray}
Here $(\tilde{a}_{kl},\tilde{b}_k)$ is a pair of Butcher's tableaux
\cite{Butcher}, of diagonally explicit Runge-Kutta schemes. This
finally completes the description of the deterministic numerical
method.
\begin{remark}
 In the case when  $\beta$ is Lipschitz but possibly degenerate,  the authors of \cite{17}, showed the $L^1$-convergence of a semi-discrete
in time relaxed scheme, see Theorem 1, Section 3 in \cite{17}. In
fact, they extended the proof of \cite{Berger_Brezis_Rogers} to the
case of a $\nu$-stages Runge-Kutta scheme. Moreover, \cite{17}
provided the following stability condition of parabolic type,
\begin{equation}\label{CFL_Cond}
    \Delta t \leq C {\Delta x}^2,
\end{equation}
where, C is a constant depending on $\beta$. At the best of our
knowledge, no such results are available in the case where $\beta$
is not Lipschitz.
\end{remark}
\section{Numerical experiments} \label{Section: Numerical
Experiments}  We  use a Matlab implementation to simulate both the
deterministic and probabilistic solutions. Concerning the plug-in
bandwidth selection procedure described in Section \ref{Section:Num.
Impl.} and based on \cite{SJ}, we have improved the code produced by
J. S.  Marron and available on
\url{http://www.stat.unc.edu/faculty/marron/marron_software.html},
by speeding up the root-finding algorithm used to solve
\eqref{AMISE-OptEps}. Furthermore, the deterministic numerical
solutions are performed using the ENO spatial reconstruction of
order 3 and a third order explicit Runge-Kutta scheme for time
stepping. We point out that the deterministic time step, denoted
from now on by $\Delta t_{det}$, is chosen with respect to the
stability condition \eqref{CFL_Cond}.
\subsection{The Classical porous media equation}
We recall that when $\beta(u)=u.|u|^{m-1},~m>1$, the PDE in
(\ref{EDP}) is nothing else but the classical porous media equation
(PME). The first numerical experiments discussed here, will be for
the mentioned $\beta$. Indeed, in the case when the initial
condition $u_0$ is a delta Dirac function at zero, we have an exact
solution provided in \cite{16}, known as the  {\it density of
Barenblatt-Pattle} and given by the following explicit formula,
\begin{eqnarray}\label{Barenblatt}
      \begin{aligned}
U(t,x)=t^{-\beta}(C-\kappa x^2 t^{-2 \beta})^{\frac{1}{m-1}}_+,
~~~x\in \mathbb{R},\ t>0,
      \end{aligned}
\end{eqnarray}
where
\[\beta=\frac{\displaystyle{1}}{\displaystyle{m+1}},~~~\kappa=\frac{\displaystyle{m-1}}{\displaystyle{2(m+1)m}},~~~
C=\left(\frac{\displaystyle{\sqrt{\kappa}}}{\displaystyle{\gamma_m}}\right)^{\frac{2(m-1)}{m+1}},~~~\gamma_m=
\int_{-\frac{\pi}{2}}^{\frac{\pi}{2}}
\left[\cos(x)\right]^{\frac{m+1}{m-1}}.\] We would now compare the
exact solution  \eqref{Barenblatt} to an approximated probabilistic
solution. However, up to now, we are not able to perform an
efficient bandwidth selection procedure in the case when the initial
condition is the law of a deterministic random variable. Since we
are nevertheless interested in exploiting \eqref{Barenblatt}, we
considered a time translation of the exact solution $U$ defined as
follows
\begin{equation}\label{BarenblattTranslate}
v(t,x)=U(t+1,x) \ \ \ \ \ \forall x\in \mathbb{R},\ \ \forall t\geq
0.
\end{equation}
Note that one can immediately deduce from
\eqref{BarenblattTranslate}, that $v$  still solves the PME but now
with a smooth initial condition given by
\begin{equation}\label{BarenblattTranslateInit}
v_0(x)=U(1,x) \ \ \ \ \ \forall x\in \mathbb{R}.
\end{equation}
In fact, in the case when the exponent $m$ is equal to $3$, the
exact solution $v$ of the PME with initial condition $v_0(x)=U(1,x)$
is given by the following explicit formula,
\begin{eqnarray}\label{BarenblattTrans_m=3}
v(t,x) = \begin{cases}
       (t+1)^{-\frac{{1}}{{4}}}\sqrt{\frac{\displaystyle{1}}{\displaystyle{\pi\sqrt{3}}}-\frac{\displaystyle{x^2}}{\displaystyle{12\sqrt{t+1}}} \  }& \text{ if }
       |x|\leq(t+1)^{\frac{{1}}{{4}}}\sqrt{\frac{\displaystyle{2}}{\displaystyle{\pi}}},\\ \\
       0 & \text{otherwise}.
       \end{cases}
\end{eqnarray}

\begin{figure}[h!]
\begin{minipage}{16cm}
 \hspace*{-0.5in}\includegraphics[width=15cm,height=12cm]{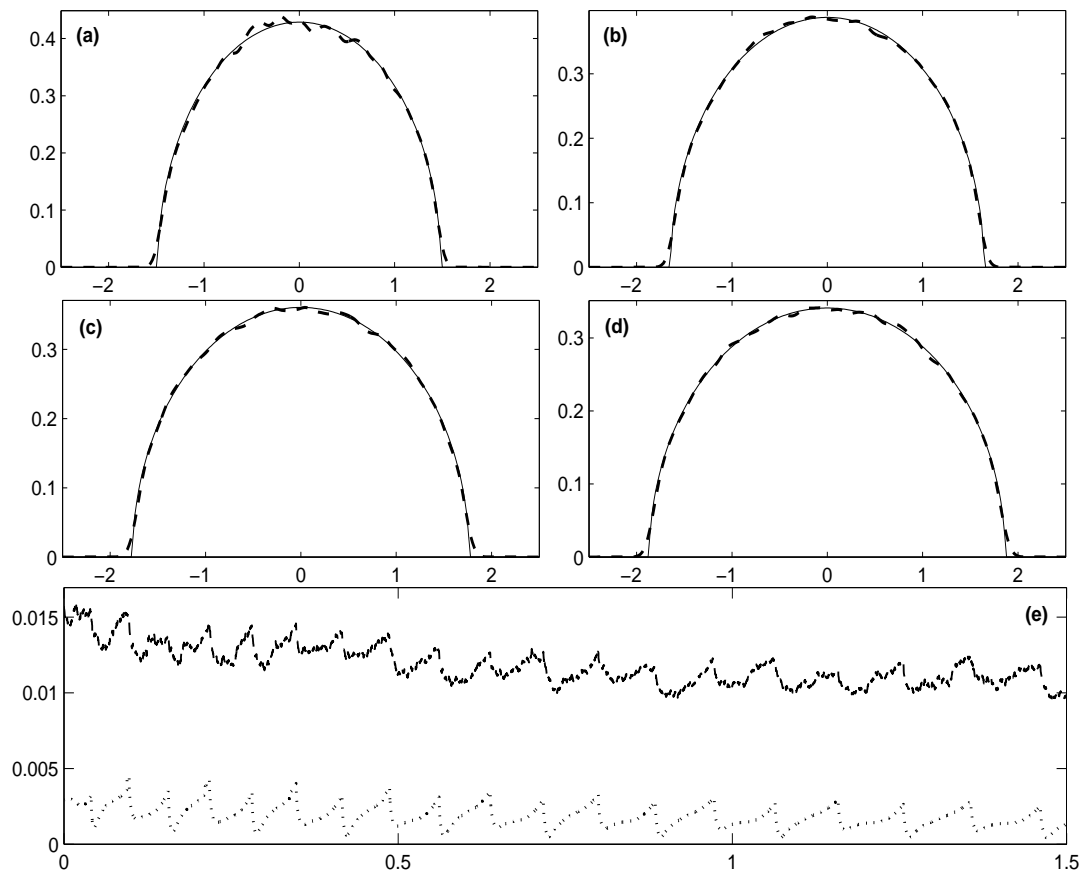}
 \end{minipage}
\begin{minipage}{12cm}
\hspace*{-0.5in} \hspace*{-0.5in}\caption{\footnotesize \textbf{-
Deterministic (doted line), probabilistic (dashed line) and exact
solutions (solid line) values at t=0 (a),
   t=0.5 (b), t=1 (c) and t=1.5 (d). The evolution of the $L^2$ deterministic
 (dote line) and probabilistic (dashed line) errors
    over the time interval $[0,1.5]$
   (e).}}
\label{fig:PME_DET_PROB_Sol_Err}
\end{minipage}
\end{figure}

 \textbf{Simulation
experiments:} we first compute both the deterministic and
probabilistic numerical solutions over the time-space grid
 $[0,1.5]\times[-2.5,2.5]$, with space step $\Delta x= 0.02$.
 We set $\Delta t_{det}=4 \times
    10^{-6}$, while, we use  $n=50000$ particles and a time step
    $\Delta t=2\times 10^{-4}$, for the
    probabilistic simulation. Figures  \ref{fig:PME_DET_PROB_Sol_Err}.(a)-(b)-(c)-(d), display the
exact and the numerical  (deterministic and probabilistic) solutions
at times $t=0$, $t=0.5$, $t=1$ and $t=T=1.5$ respectively. The exact
solution of the PME, defined in (\ref{BarenblattTrans_m=3}), is
depicted by solid lines.

Besides, Figure \ref{fig:PME_DET_PROB_Sol_Err}.(e) describes the
time evolution of both the discrete $L^2$ deterministic and
probabilistic errors  on the time interval $[0,1.5]$.

The $L^1$ errors behave very similarly as well in the present case
as in the Heaviside case, treated in subsection
\ref{Subsec:Heaviside Case}.

\subsection{The Heaviside case} \label{Subsec:Heaviside Case}
The second family of numerical experiments discussed here, concerns
 $\beta$ defined by \eqref{Heav}. Since we do not have an exact
solution of the diffusion problem \eqref{EDP}, for the mentioned
$\beta$, we decided to compare the probabilistic solution to the
approximation obtained via the deterministic algorithm described in
Section \ref{Section:Det. App}. Indeed, we shall simulate both
solutions according to several types of initial data $u_0$ and with
different values of the critical threshold $u_c$.

Empirically, after various experiments, it appears that for a fixed
threshold $u_c$, the numerical solution approaches some limit
function which seems to belong to the "attracting" set
\begin{equation}\label{AttractSet}
    \mathcal{J}=\{f \in L^1(\mathbb{R})| \int f(x)dx=1,\ \ |f|\leq
    u_c\};
\end{equation}
in fact $\mathcal{J}$ is the closure in $L^1$ of $\mathcal{J}_0=\{f
: \mathbb{R} \rightarrow \mathbb{R}_+ | \  \beta(f)=0\}$. At this
point,  the following theoretical questions arise.
\begin{description}
    \item [(1)]Does indeed $u(t,\cdot)$ have a limit $u_{\infty}$ when $t \rightarrow
    \infty$?
    \item [(2)]If yes does $u_{\infty}$  belong to $\mathcal{J}$?
    \item [(3)]If (2) holds, do we have $u(t,\cdot)=u_{\infty}$
    for $t$ larger than a finite time $\tau$?
\end{description}
A similar behavior was observed for different $\beta$ which are
strictly increasing after some zero.

\subsubsection{Trimodal initial condition}
For the $\beta$   given by \eqref{Heav}, we consider  an initial
condition being a mixture of three Gaussian densities with three
modes at some distance from each other, i.e.
\begin{equation}\label{CondInitTrino}
u_0(x)=\frac{\displaystyle{1}}{\displaystyle{3}}\left(p(x,\mu_1,\sigma_1)+p(x,\mu_2,\sigma_2)+p(x,\mu_3,\sigma_3)\right),
\end{equation}
where,
\begin{equation}\label{Gaussian}
p(x,\mu,\sigma)=\frac{\displaystyle{1}}{\displaystyle{\sqrt{2\pi}\sigma}}\exp(-\frac{\displaystyle{(x-\mu)^2}}{\displaystyle{2\sigma^2}}).
\end{equation}

\textbf{Simulation experiments:} for this specific type of initial
condition $u_0$, we consider two test cases depending on the value
taken by the critical threshold  $u_c$. We set, for instance,
$\mu_1=-\mu_3=-4$, $\mu_2=0$ and $\sigma_1=0.1$, $\sigma_2=0.2$,
$\sigma_3=0.3$. \\

\begin{figure}[p]
\begin{minipage}{16cm}
\hspace*{-0.5in}\includegraphics[width=16cm,height=10cm]{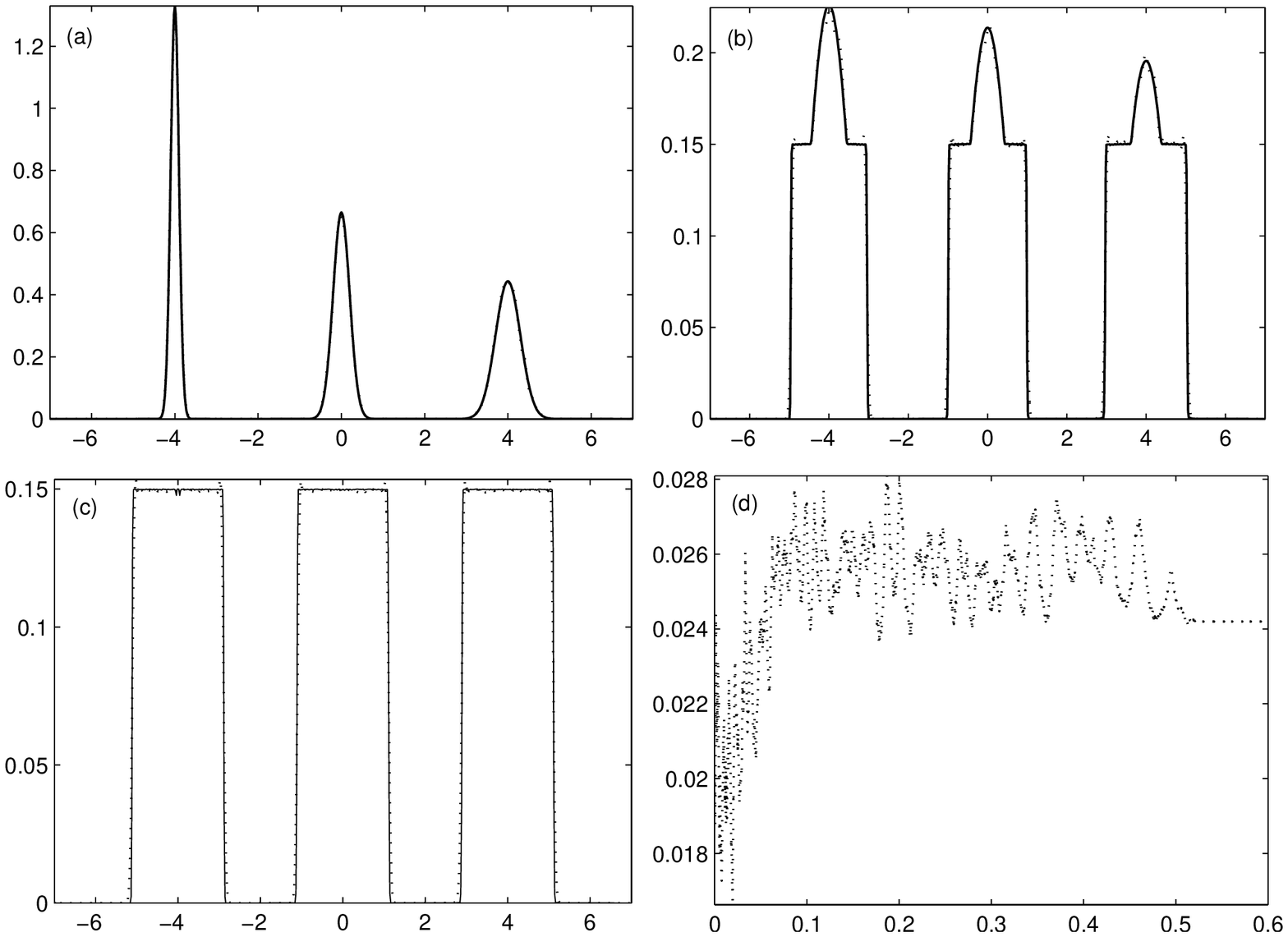}
\end{minipage}
%%%%%%%%%%%
\vspace{0.01in}
%%%%%%%%%%%
\begin{minipage}{14cm}
\hspace*{-0.5in}\caption{\footnotesize \textbf{- Test case 1:
Deterministic (solid line) and probabilistic (doted line) solution
  values at t=0 (a),
   t=0.3 (b), t=0.6 (c). The evolution of the $L^2$-norm of the difference  over the time interval $[0,0.6]$
   (d).}}
\label{fig:TrinoToTrino_Sol_Err}
\end{minipage}
%%%%%%%%%%%
\vspace{0.05in}
%%%%%%%%%%%
\begin{minipage}{16cm}
\hspace*{-0.5in}\includegraphics[width=16cm,height=10cm]{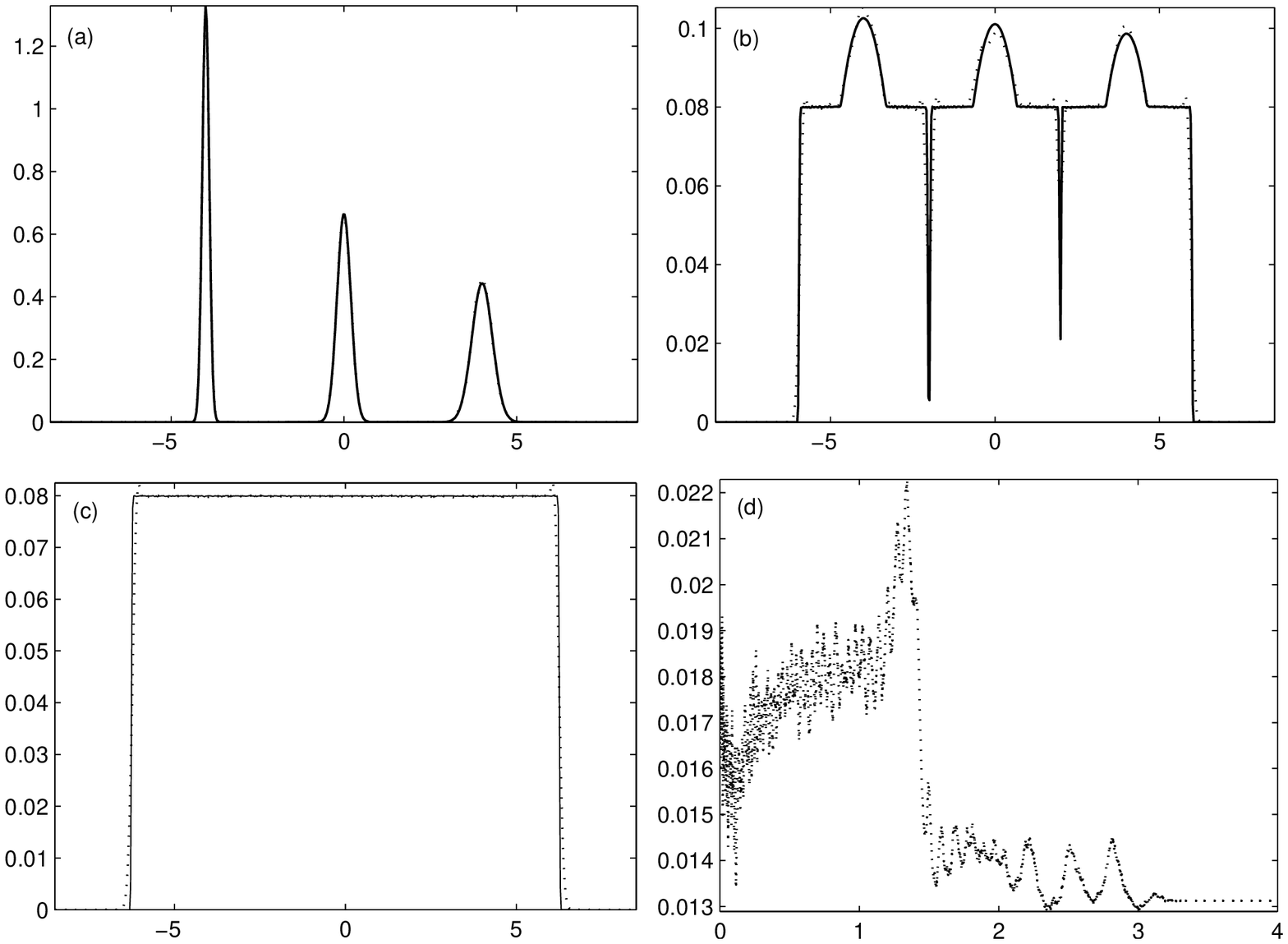}
\end{minipage}
%%%%%%%%%%%
\vspace{0.01in}
%%%%%%%%%%%
\begin{minipage}{14cm}
\hspace*{-0.5in}\caption{\footnotesize \textbf{- Test case 2:
Deterministic (solid line) and probabilistic (doted line) solution
  values at t=0 (a),
   t=2 (b), t=4 (c). The evolution of the $L^2$-norm of the difference  over the time interval $[0,4]$
   (d).}}
\label{fig:TrinoToUni_Sol_Err}
\end{minipage}
\end{figure}

\textbf{Test case 1} :  we start with $u_c=0.15$, and
    a time-space grid $[0,0.6]\times[-7,7]$, with a space step $\Delta x=0.02$. For the deterministic approximation,
    we set $\Delta t_{det}=4 \times
    10^{-6}$.  The probabilistic simulation  uses
     $n=50000$ particles and a time step $\Delta t=2\times 10^{-4}$.  Figures \ref{fig:TrinoToTrino_Sol_Err}.(a)-(b)-(c), displays both
    the deterministic and probabilistic numerical solutions at times $t=0$, $t=0.3$ and $t=T=0.6$, respectively. On the other hand, the time evolution of
     the   $L^2$-norm of the difference between the two numerical solutions is depicted in Figure
     \ref{fig:TrinoToTrino_Sol_Err}.(d).\\

\textbf{Test case 2} : we choose now as critical value $u_c=0.08$
and a time-space grid $[0,4]\times[-8.5,8.5]$, with a space step
$\Delta x=0.02$. We set  $\Delta t_{det}=4\times 10^{-6}$
    and the probabilistic approximation is performed using $n=50000$ particles and a time step
    $\Delta t=2\times 10^{-4}$.  Figures \ref{fig:TrinoToUni_Sol_Err}.(a)-(b)-(c) and  \ref{fig:TrinoToUni_Sol_Err}.(d), show respectively
    the numerical (probabilistic and deterministic) solutions and the  $L^2$-norm of the difference between
    the two.
\subsubsection{Uniform and Normal densities mixture initial condition}
We proceed with  $\beta$ given by \eqref{Heav}. We are now
interested in an initial condition $u_0$, being a mixture of a
Normal and an Uniform density, i.e.,
\begin{equation}\label{CondInitMixtureUnifNorm}
u_0(x)=\frac{\displaystyle{1}}{\displaystyle{2}}\left(p(x,-1,0.2)+\mathds{1}_{[0,1]}(x)\right),
\end{equation}
where, $p$ is defined in \eqref{Gaussian}.\\

\textbf{Simulation experiments:}\\

\textbf{Test case 3} : we perform both the approximated
deterministic and
    probabilistic solutions in the case where $u_c=0.3$, on
    the time-space grid $[0,0.5]\times[-2.5,2]$, with a space step $\Delta x=0.02$. We use  $n=50000$ particles and a time step
    $\Delta t=2\times 10^{-4}$, for the
    probabilistic simulation. Moreover, we set $\Delta t_{det}=4 \times
    10^{-6}$. Figures
    \ref{fig:Mixte_Unif_Norm_Sol_Err}.(a)-(b)-(c) illustrate
    those
    approximated solutions at times $t=0$, $t=0.1$ and  $t=T=0.5$. Furthermore, we compute the  $L^2$-norm of the
difference between the numerical deterministic solution and the
probabilistic one. Values of this error, are displayed in Figure
\ref{fig:Mixte_Unif_Norm_Sol_Err}.(d), at each probabilistic time
step.

\subsubsection{Uniform densities  mixture initial condition}
Now, with $\beta$ given by \eqref{Heav}, we consider an initial
condition $u_0$ being a mixture of Uniform densities, i.e.,
\begin{equation}\label{CondInitSomUnif}
u_0(x)=\frac{\displaystyle{1}}{\displaystyle{5}}\mathds{1}_{[0,1]}(x)+
\frac{\displaystyle{3}}{\displaystyle{4}}\mathds{1}_{[-\frac{{1}}{{5}},\frac{{1}}{{5}}]}(x)
+\frac{\displaystyle{5}}{\displaystyle{8}}\mathds{1}_{[\frac{{6}}{{5}},2]}(x),
\end{equation}

\textbf{Simulation experiments:}\\

\begin{figure}[p]
\begin{minipage}{16cm}
\hspace*{-0.5in}\includegraphics[width=16cm,height=10cm]{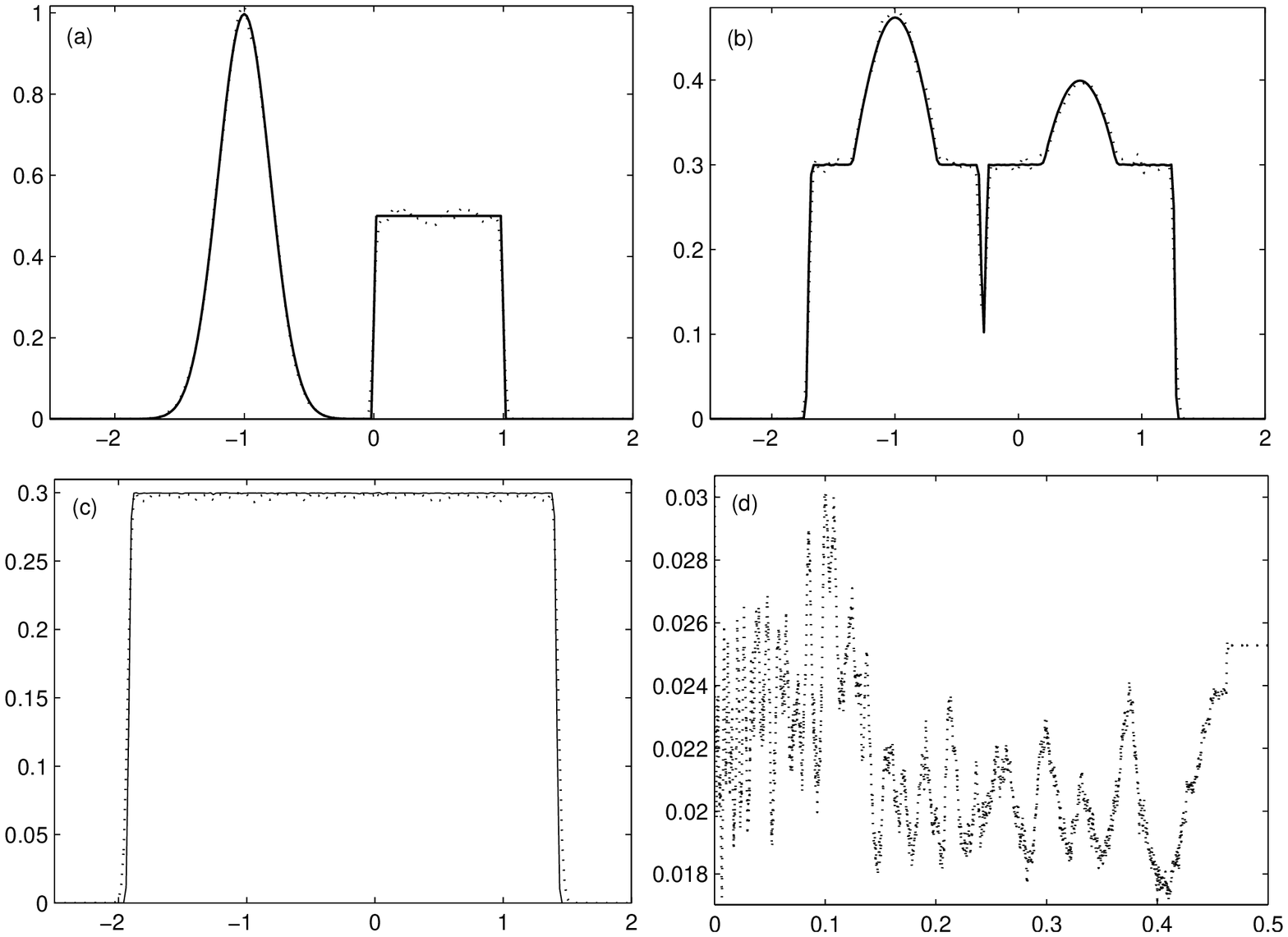}
\end{minipage}
%%%%%%%%%%%
\vspace{0.01in}
%%%%%%%%%%%
\begin{minipage}{14cm}
\hspace*{-0.5in}\caption{\footnotesize \textbf{- Test case 3:
Deterministic (solid line) and probabilistic (doted line) solutions
  values at t=0 (a),
   t=0.1 (b), t=0.5 (c). The evolution of the $L^2$-norm of the difference  over the time interval $[0,0.5]$
   (d).}}
\label{fig:Mixte_Unif_Norm_Sol_Err}
\end{minipage}
%%%%%%%%%%%
\vspace{0.05in}
%%%%%%%%%%%
\begin{minipage}{16cm}
\hspace*{-0.5in}\includegraphics[width=16cm,height=10cm]{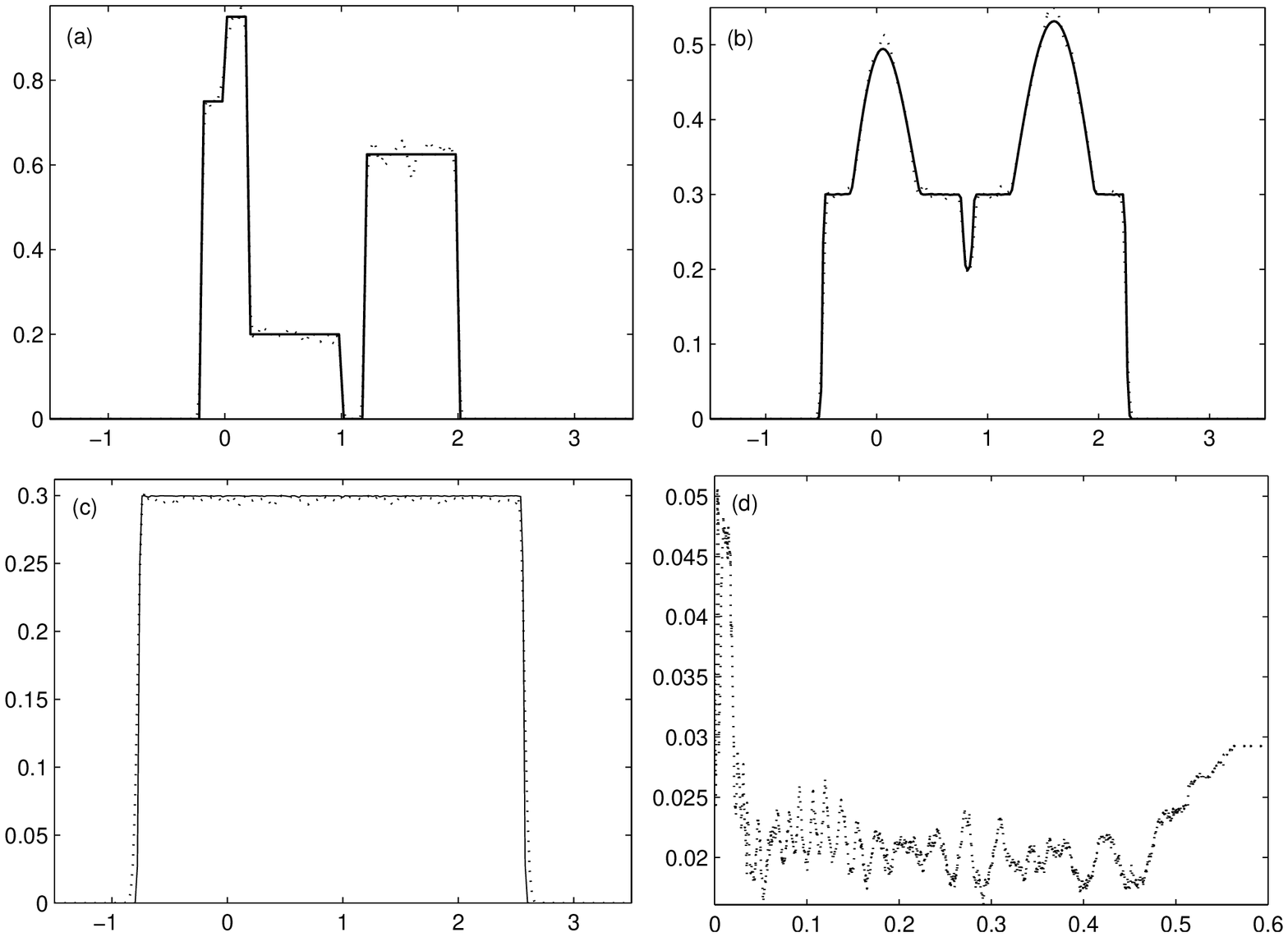}
\end{minipage}
%%%%%%%%%%%
\vspace{0.01in}
%%%%%%%%%%%
\begin{minipage}{14cm}
\hspace*{-0.5in}\caption{\footnotesize \textbf{- Test case 4:
Deterministic (solid line) and probabilistic solution
 (doted line) values at t=0 (a),
   t=0.1 (b), t=0.6 (c). The evolution of the $L^2$-norm of the difference  over the time interval $[0,0.6]$
   (d).}}
\label{fig:Somme_Unif_Sol_Err}
\end{minipage}
\end{figure}

 \textbf{Test case 4} : we approximate  the  deterministic and
    probabilistic solutions in the case where $u_c=0.3$, on
    the time-space grid $[0,0.6]\times[-1.5,3.5]$, with a space step $\Delta x=0.02$.
   The deterministic time step $\Delta t_{det}=4\times 10^{-6}$,
   while
    the probabilistic solution is computed using $n=50000$ particles and a time step
    $\Delta t=2\times 10^{-4}$. We illustrate in Figure \ref{fig:Somme_Unif_Sol_Err}.(a)-(b)-(c), both the deterministic and probabilistic
    numerical solutions at  times  $t=0$, $t=0.1$ and $t=T=0.6$ respectively; while the time evolution of the
     $L^2$-norm of the difference between them, is shown
in Figure \ref{fig:Somme_Unif_Sol_Err}.(d).

\subsubsection{Square root initial condition} Finally, the last test
case concerns an initial condition $u_0$ defined as follows :
\begin{equation}\label{CondInitRacineCarre}
u_0(x)=\frac{\displaystyle{3}}{\displaystyle{4}}\sqrt{|x|}\mathds{1}_{[-1,1]}(x).
\end{equation}

\textbf{Simulation experiments:}\\

 \textbf{Test case 5} : we simulate the  probabilistic and deterministic
 solutions over the time-space grid $[0,0.45]\times[-2,2]$, using a
 space step $\Delta x=0.02$ and setting the
 critical threshold $u_c=0.35$. Moreover, the deterministic time
 step $\Delta t_{det}=4\times 10^{-6}$. On the other hand, we use $n=50000$ particles and a time step
    $\Delta t=2\times 10^{-4}$ to compute the probabilistic
    approximation.

Figures \ref{fig:Puiss_Gamma_Sol_Err}.(a)-(b)-(c), show both the
    deterministic and probabilistic numerical solutions at times $t=0$, $t=0.04$ and
    $t=T=0.45$, respectively.

    The evolution of the $L^2$-norm of the difference between
    these two solutions, over the time interval $[0,0.45]$, is depicted in Figure
    \ref{fig:Puiss_Gamma_Sol_Err}.(d).

\subsection{Concluding remarks}
\begin{figure}[p]
\begin{minipage}{16cm}
\hspace*{-0.5in}\includegraphics[width=16cm,height=10cm]{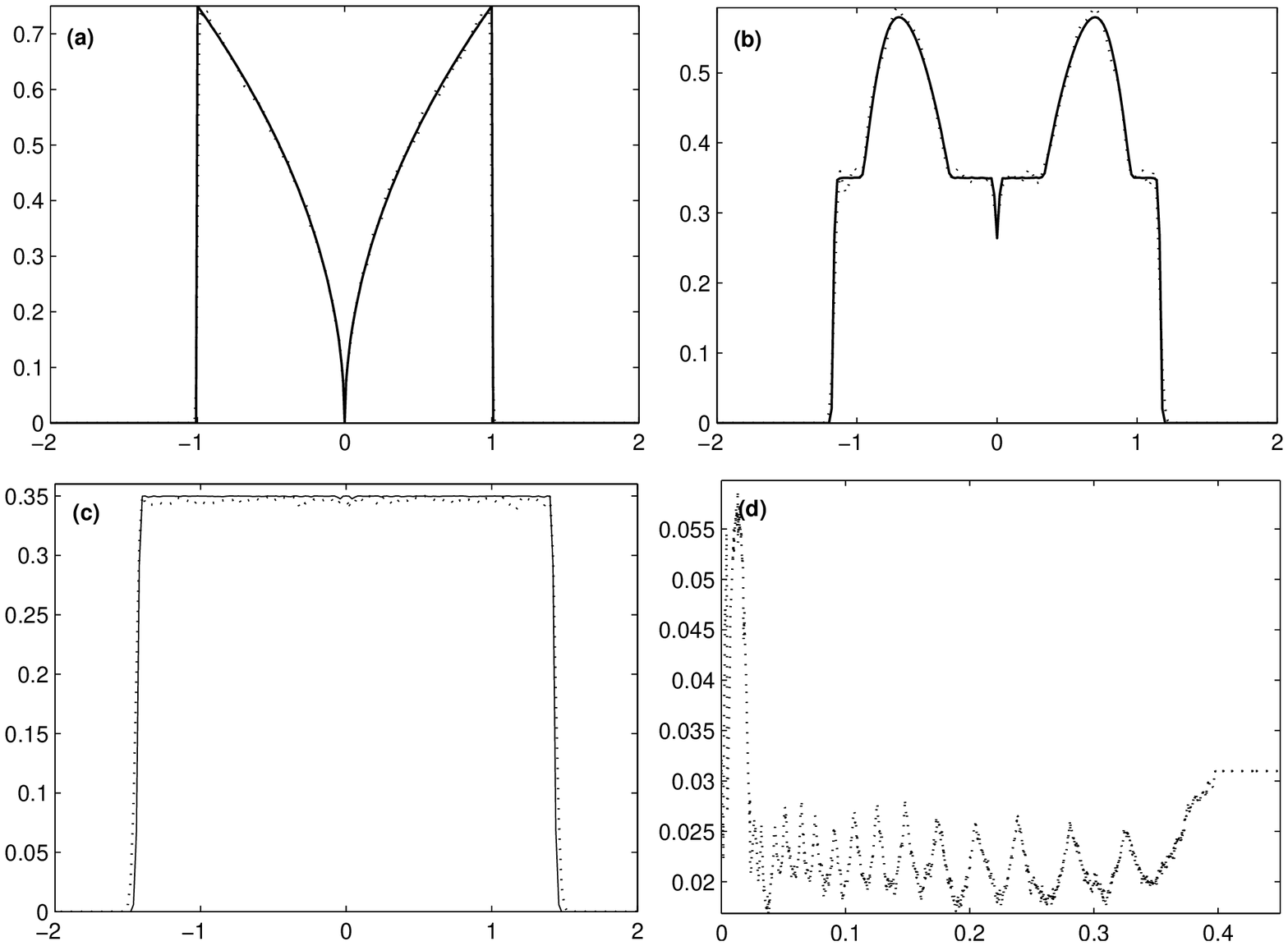}
\end{minipage}
%%%%%%%%%%%
\vspace{0.01in}
%%%%%%%%%%%
\begin{minipage}{14cm}
\hspace*{-0.5in}\caption{\footnotesize \textbf{- Test case 5:
Deterministic (solid line) and probabilistic (doted line) solution
  values at t=0 (a),
   t=0.04 (b), t=0.45 (c). The evolution of the $L^2$-norm of the difference  over the time interval $[0,0.45]$
   (d).}}
\label{fig:Puiss_Gamma_Sol_Err}
\end{minipage}
%%%%%%%%%%%
\vspace{0.05in}
%%%%%%%%%%%
\begin{minipage}{16cm}
\hspace*{-0.4in}\includegraphics[width=16cm,height=10cm]{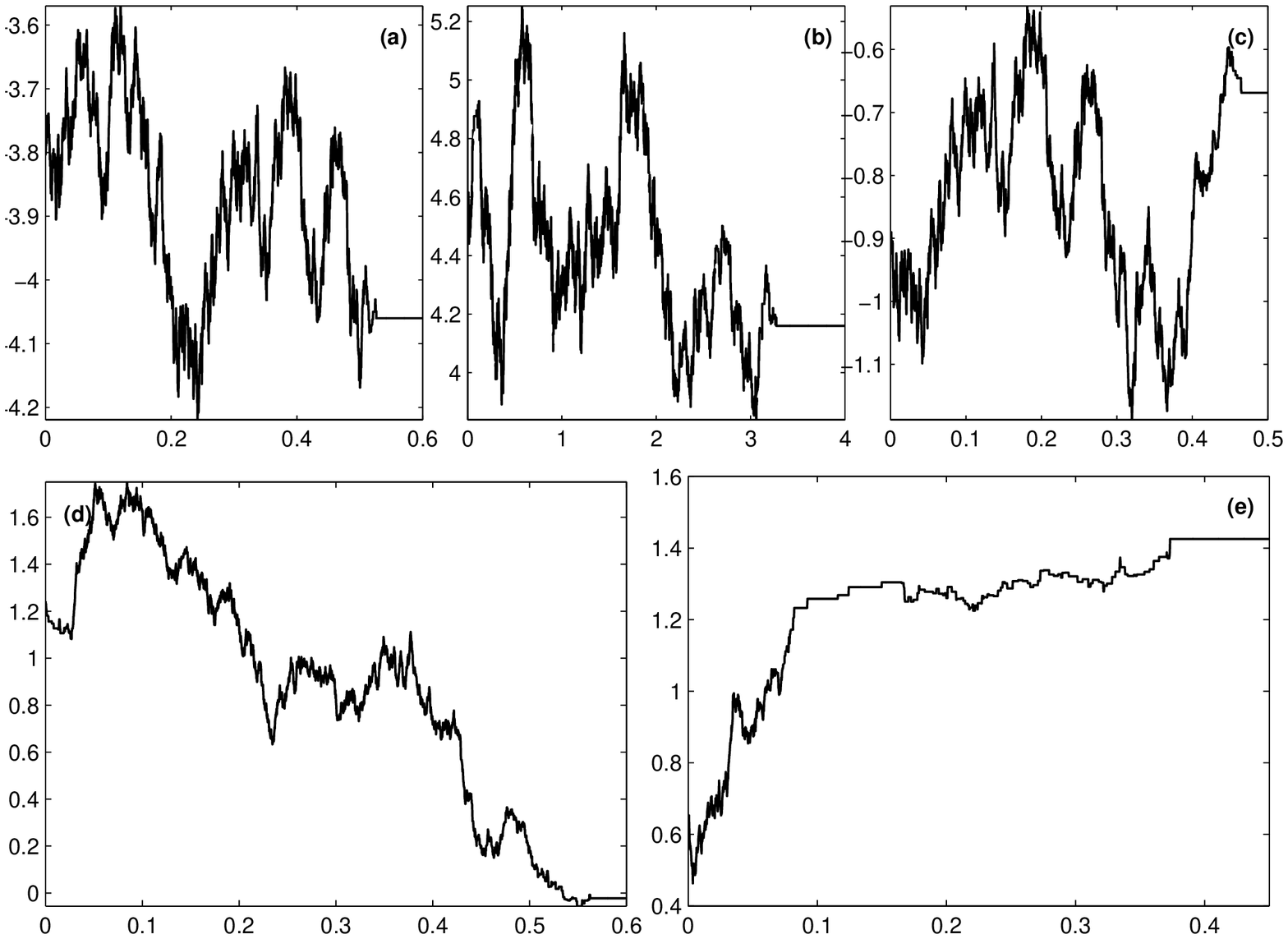}
\end{minipage}
%%%%%%%%%%%
\vspace{0.01in}
%%%%%%%%%%%
\begin{minipage}{14cm}
\hspace*{-0.5in}\caption{\footnotesize \textbf{- Representation of a
process trajectory
 for the Test case 1 (a), Test case 2 (b), Test case 3
(c), Test case 4 (d) and  Test case 5 (e),  respectively.}}
\label{fig:Trajectories}
\end{minipage}
\end{figure}

\begin{description}
    \item [(1)] Figure \ref{fig:Trajectories}, displays a single trajectory for
each one of the test cases described above. In fact, we observe
that, in all cases, the process trajectory stops not later than the
instant of stabilization of the macroscopic
distribution.\\
    \item [(2)] We have performed deterministic and probabilistic
 numerical solutions for \eqref{EDP}, with a coefficient $\beta$ defined by \eqref{Heav}. Even though
 the procedures being used were different, the simulation experiments clearly
 show that both methods produce very close approximated  solutions,  all over the considered time
 interval.\\
    \item [(3)] We point out that, the error committed  by the Monte
    Carlo simulations largely dominates the one related to the Euler
    scheme. Consequently, the choice of the probabilistic time step is not so
    important.\\
    \item [(4)] The probabilistic algorithm can be
    parallelized on a Graphical Processor Unit (GPU), such that we can speed-up its time machine execution; on the other hand, for the deterministic algorithm,
    this transformation is far from being obvious, see \cite{TR_Cuvelier}.

\end{description}

%---------------------------------------------------------------------------------------------------------------%
%***************************************************************************************************************%
%---------------------------------------------------------------------------------------------------------------%
%                                                                                                               %
%         ANNEXES ANNEXES ANNEXES ANNEXES ANNEXES ANNEXES ANNEXES ANNEXES ANNEXES ANNEXES ANNEXES ANNEXES       %
%                                                                                                               %
%---------------------------------------------------------------------------------------------------------------%
%***************************************************************************************************************%
%---------------------------------------------------------------------------------------------------------------%

\newpage
\newpage
\section{Annexes}\label{Annex}
 Let
$\mathbf{V}\in\mathbb{R}^{N_x}$ such that $V_i=v(x_i)$,  $\forall
i\in\left[\hspace{-1ex} \left[\hspace{0.5ex}1,N_x\hspace{0.5ex}
\right]\hspace{-1ex}\right]$, where $v$ is a function defined on
$[a,b]$. Note that the points $(x_i)$ are still defined as in
Section \ref{Section:Det. App}. Moreover,
$\mathcal{M}_{m,n}(\mathbb{R})$ denotes the linear space  of real
matrices with $m$ rows and $n$ columns.
\subsection{Interpolating polynomial of a function\label{SecAnnex1}}
We aim to approximate  $ v(x_{i+1/2})$ and $ v(x_{i-1/2})$ for every
$i\in \left[\hspace{-1ex} \left[\hspace{0.5ex}1,N_x\hspace{0.5ex}
\right]\hspace{-1ex}\right]$ . In order to do this, we use properly
chosen Lagrange interpolation polynomials of degree $k-1$ .

On every interval (or cell)  $I_i=[x_{i-1/2},x_{i+1/2}]$, with $i
\in \left[\hspace{-1ex} \left[\hspace{0.5ex}1,N_x\hspace{0.5ex}
\right]\hspace{-1ex}\right]$ ,  we construct an interpolation
polynomial $\mathbb{P}_{k,i}$ by selecting $k$ consecutive points
containing $x_i$ : the so-called \emph{stencil} denoted by
\[S(i)\equiv\{I_{i-r},\ldots,I_{i+s}\}\]  and defined by
$\{x_{i-r},x_{i-r+1},\ldots,x_{i+s-1},x_{i+s}\}$,  where $r,s$ are
positive integers and $r+s+1=k$.  We denote by $R(i)$, the value
taken by $r$ for the interval $I_i$ with
 an ENO stencil, see \cite{Shu}.

The Lagrange interpolation polynomial of degree $k-1$, on the
interval $I_i$, associated to the stencil $S(i)$ is then given by :
\begin{equation}\label{Annex(4.1)}
\mathbb{P}_{k,i}^{[r]}(x)=\sum_{j=0}^{k-1} V_{i-r+j} L_j^{[r]}(x),\
\ \ \ \forall x\in I_i,
\end{equation}
where,
\begin{equation}\label{Lagrange_Poly}
L_j^{[r]}(x)=\prod_{ \substack{l=0\\l\neq
j}}^{k-1}\frac{x-x_{i-r+l}}{x_{i-r+j}-x_{i-r+l}}.
\end{equation}
Now, we need to compute the  polynomial defined in
\eqref{Lagrange_Poly} at the points $x_{i-1/2}$ and $x_{i+1/2}$. In
fact, since the points are equidistant, we have for every $ (r,j)\in
{\left[\hspace{-1ex}\left[\hspace{0.5ex}0,{k-1}\hspace{0.5ex}
\right]\hspace{-1ex}\right]}^2$ ,

\begin{equation*}\label{Lagrange_Poly_Bis}
\left\{
 \begin{array}{rcl}
 L_j^{[r]}(x_{i-1/2})&=&\ \prod\limits_{ \substack{l=0\\ l\neq
j}}^{k-1} \frac{\displaystyle{r-l-1/2}}{\displaystyle{j-l}},\\ \\
L_j^{[r]}(x_{i+1/2})&=&\ \prod\limits_{ \substack{l=0\\ l\neq
j}}^{k-1} \frac{\displaystyle{r-l+1/2}}{\displaystyle{j-l}}.
\end{array}
\right.
\end{equation*}
Then, we define $\mathbb{C}\in \mathcal{M}_{k+1,k}(\mathbb{R})$,  as
follows
\begin{equation}\label{Annex(4.2)}
\mathbb{C}_{r+1,j+1}=\prod_{ \substack{l=0\\ l\neq j}}^{k-1}
\frac{r-l-1/2}{j-l},\ \ \ \forall (r,j)\in
\left[\hspace{-1ex}\left[\hspace{0.5ex}0,{k}\hspace{0.5ex} \right]
\hspace{-1ex}\right]\times\left[\hspace{-1ex}\left[\hspace{0.5ex}0,{k-1}\hspace{0.5ex}
\right]\hspace{-1ex}\right].
\end{equation}
Substituting \eqref{Annex(4.2)} in \eqref{Annex(4.1)} and using the
ENO stencil, we get  $\forall i \in \left[\hspace{-1ex}
\left[\hspace{0.5ex}1,N_x\hspace{0.5ex}
\right]\hspace{-1ex}\right]$,
\begin{eqnarray}
v(x_{i-1/2})\approx v_{i-1/2}^+=
\mathbb{P}_{k,i}^{[R(i)]}(x_{i-1/2})= \sum_{j=0}^{k-1} V_{i-R(i)+j}
\mathbb{C}_{R(i)+1,j+1},\label{InterpData1}\\
v(x_{i+1/2})\approx
v_{i+1/2}^-=\mathbb{P}_{k,i}^{[R(i)]}(x_{i+1/2})= \sum_{j=0}^{k-1}
V_{i-R(i)+j} \mathbb{C}_{R(i)+2,j+1}.\label{InterpData2}
\end{eqnarray}
\subsection{Interpolation polynomial for the derivative of a function \label{SecAnnex2}}
Now, we would like  to approximate
$\frac{\displaystyle{dv}}{\displaystyle{dx}}(x_i)$,
$\frac{\displaystyle{dv}}{\displaystyle{dx}}(x_{i-1/2})$ and
$\frac{\displaystyle{dv}}{\displaystyle{dx}}(x_{i+1/2})$, for every
$ i \in\left[\hspace{-1ex}\left[\hspace{0.5ex}1,{N_x}\hspace{0.5ex}
\right]\hspace{-1ex}\right] $ . In fact, deriving equation
(\ref{Annex(4.1)}), implies
\begin{equation}\label{Annex(4.5)}
\frac{\displaystyle{d\mathbb{P}_{k,i}^{[r]}}}{\displaystyle{dx}}(x)=\sum_{j=0}^{k-1}
V_{i-r+j} \frac{\displaystyle{dL_j^{[r]}}}{\displaystyle{dx}}(x),\ \
\ \ \forall x\in I_i.
\end{equation}
On the other hand,  for every $j\in
\left[\hspace{-1ex}\left[\hspace{0.5ex}0,{l-1}\hspace{0.5ex}
\right]\hspace{-1ex}\right]$ , we have
\begin{equation}\label{LjDstencil}
\frac{\displaystyle{dL_j^{[r]}}}{\displaystyle{dx}}(x)=\frac{\displaystyle
  \sum_{ \substack{m=0\\ m\neq j}}^{k-1}
  \prod_{ \substack{l=0\\ l\neq j,m}}^{k-1} (x-x_{i-r+l})}
{\displaystyle  \prod_{ \substack{l=0\\ l\neq j}}^{k-1}
(x_{i-r+j}-x_{i-r+l})},\ \ \forall x\in I_i.
\end{equation}
Since the points are equidistant, we get
\begin{equation*}
 \frac{\displaystyle{dL_j^{[r]}}}{\displaystyle{dx}}(x_i)=\
\frac{\displaystyle
  \sum\limits_{ \substack{m=0\\ m\neq j}}^{k-1}
  \prod\limits_{ \substack{l=0\\ l\neq j,m}}^{k-1} (r-l)}
{\displaystyle  \Delta x\prod\limits_{ \substack{l=0\\ l\neq
j}}^{k-1} (j-l)}
\end{equation*}
and
\begin{equation*}
\frac{\displaystyle{dL_j^{[r]}}}{\displaystyle{dx}}(x_{i-1/2})=\
\frac{\displaystyle
  \sum\limits_{ \substack{m=0\\ m\neq j}}^{k-1}
  \prod\limits_{ \substack{l=0\\ l\neq j,m}}^{k-1} (r-l-1/2)}
{\displaystyle  \Delta x\prod\limits_{ \substack{l=0\\ l\neq
j}}^{k-1} (j-l)},\ \
\frac{\displaystyle{dL_j^{[r]}}}{\displaystyle{dx}}(x_{i+1/2})= \
\frac{\displaystyle
  \sum\limits_{ \substack{m=0\\ m\neq j}}^{k-1}
  \prod\limits_{ \substack{l=0\\ l\neq j,m}}^{k-1} (r-l+1/2)}
{\displaystyle  \Delta x\prod_{ \substack{l=0\\ l\neq j}}^{k-1}
(j-l)}.
\end{equation*}
Then, we define $\mathbb{D}\in \mathcal{M}_{k,k}(\mathbb{R})$  by
\begin{equation}\label{D_matrice}
\mathbb{D}_{r+1,j+1}=\frac{\displaystyle
  \sum_{ \substack{m=0\\ m\neq j}}^{k-1}
  \prod_{ \substack{l=0\\ l\neq j,m}}^{k-1} (r-l)}
{\displaystyle  \Delta x\prod_{ \substack{l=0\\ l\neq j}}^{k-1}
(j-l)},\ \ \ \ \ \forall (r,j) \in
\left[\hspace{-1ex}\left[\hspace{0.5ex}0,{k-1}\hspace{0.5ex}
\right]\hspace{-1ex}\right]^2,
\end{equation}
and $\overline{\mathbb{D}}\in \mathcal{M}_{k+1,k}(\mathbb{R})$  as
follows
\begin{equation}\label{D_bar_matrice}
\overline{\mathbb{D}}_{r+1,j+1}=\frac{\displaystyle
  \sum_{ \substack{m=0\\ m\neq j}}^{k-1}
  \prod_{ \substack{l=0\\ l\neq j,m}}^{k-1} (r-l-1/2)}
{\displaystyle  \Delta x\prod_{ \substack{l=0\\ l\neq j}}^{k-1}
(j-l)},\ \ \ \  \forall (r,j) \in
\left[\hspace{-1ex}\left[\hspace{0.5ex}0,{k}\hspace{0.5ex} \right]
\hspace{-1ex}\right]\times
\left[\hspace{-1ex}\left[\hspace{0.5ex}0,{k-1}\hspace{0.5ex}
\right]\hspace{-1ex}\right].
\end{equation}
Therefore, replacing  \eqref{D_matrice} and \eqref{D_bar_matrice} in
\eqref{Annex(4.5)},  for every $i \in
\left[\hspace{-1ex}\left[\hspace{0.5ex}1,N_x\hspace{0.5ex}
\right]\hspace{-1ex}\right]$ , we obtain :
\begin{eqnarray}
\frac{\displaystyle{dv}}{\displaystyle{dx}}(x_i) &\approx & dv_i=\frac{\displaystyle{d\mathbb{P}_{k,i}^{[r]}}}{\displaystyle{dx}}(x_i)= \sum_{j=0}^{k-1} V_{i-r+j} \mathbb{D}_{r+1,j+1},\label{DerivInterpData1}\\
\frac{\displaystyle{dv}}{\displaystyle{dx}}(x_{i-1/2}) &\approx & dv_{i-1/2}^+= \frac{\displaystyle{d\mathbb{P}_{k,i}^{[r]}}}{\displaystyle{dx}}(x_{i-1/2})= \sum_{j=0}^{k-1} V_{i-r+j} \overline{\mathbb{D}}_{r+1,j+1},\label{DerivInterpData2}\\
\frac{\displaystyle{dv}}{\displaystyle{dx}}(x_{i+1/2}) &\approx &
dv_{i+1/2}^-=\frac{\displaystyle{d\mathbb{P}_{k,i}^{[r]}}}{\displaystyle{dx}}(x_{i+1/2})=
\sum_{j=0}^{k-1} V_{i-r+j}
\overline{\mathbb{D}}_{r+2,j+1}.\label{DerivInterpData3}
\end{eqnarray}
\newpage
\bibliographystyle{mcma}
\bibliography{Belaribi_CR_biblio}

%% file: MainBCRSentNov2010.bbl
\providecommand{\bysame}{\leavevmode\hbox to3em{\hrulefill}\thinspace}
\providecommand{\MR}{\relax\ifhmode\unskip\space\fi MR }
% \MRhref is called by the amsart/book/proc definition of \MR.
\providecommand{\MRhref}[2]{%
  \href{http://www.ams.org/mathscinet-getitem?mr=#1}{#2}
}
\providecommand{\href}[2]{#2}
\begin{thebibliography}{10}

\bibitem{Aregba}
D.~Aregba-Driollet, R.~Natalini and S.~Tang, \emph{Explicit diffusive kinetic
  schemes for nonlinear degenerate parabolic systems}, Math. Comp. {73} (2004),
  pp.~63--94 (electronic).

\bibitem{1}
P.~Bak, \emph{How Nature Works: The science of Self-Organized Criticality}.
  Springer-Verlag New York, Inc, 1986.

\bibitem{10}
V.~Barbu, M.~R\"ockner and F.~Russo, \emph{Probabilistic representation for
  solutions of an irregular porous media type equation: the irregular
  degenerate case}. To appear: Prob. Th. Rel. Fields. Available at \url{
  http://hal.inria.fr/inria-00410248/fr/}.

\bibitem{16}
G.~I. Barenblatt, \emph{On some unsteady motions of a liquid and gas in a
  porous medium}, Akad. Nauk SSSR. Prikl. Mat. Meh. {16} (1952), pp.~67--78.

\bibitem{8}
S.~Benachour, P.~Chassaing, B.~Roynette and P.~Vallois, \emph{Processus
  associ\'es \`a\ l'\'equation des milieux poreux}, Ann. Scuola Norm. Sup. Pisa
  Cl. Sci. (4) {23} (1996), pp.~793--832 (1997).

\bibitem{4}
P.~Benilan, H.~Brezis and M.~G. Crandall, \emph{A semilinear equation in
  {$L^1(\mathbb{R}^N)$}}, Ann. Scuola Norm. Sup. Pisa Cl. Sci. (4) {2} (1975),
  pp.~523--555.

\bibitem{6}
P.~Benilan and M.~G. Crandall, \emph{The continuous dependence on {$\varphi $}
  of solutions of {$u\sb{t}-\Delta \varphi (u)=0$}}, Indiana Univ. Math. J.
  {30} (1981), pp.~161--177.

\bibitem{Berger_Brezis_Rogers}
A.~E. Berger, H.~Br{\'e}zis and J.~C.~W. Rogers, \emph{A numerical method for
  solving the problem {$u_{t}-\Delta f(u)=0$}}, RAIRO Anal. Num\'er. {13}
  (1979), pp.~297--312.

\bibitem{3}
P.~Blanchard, M.~R{\"o}ckner and F.~Russo, \emph{Probabilistic representation
  for solutions of an irregular porous media type equation}, Ann. Probab. {38}
  (2010), pp.~1870--1900.

\bibitem{BossyTalay95}
M.~Bossy and D.~Talay, \emph{A stochastic particle method for some
  one-dimensional nonlinear p.d.e}, Math. Comput. Simulation {38} (1995),
  pp.~43--50. Probabilit{\'e}s num{\'e}riques (Paris, 1992).

\bibitem{BossyTalay97}
\bysame, \emph{A stochastic particle method for the {M}c{K}ean-{V}lasov and the
  {B}urgers equation}, Math. Comp. {66} (1997), pp.~157--192.

\bibitem{BKGSyst}
F.~Bouchut, F.~R. Guarguaglini and R.~Natalini, \emph{Diffusive {BGK}
  approximations for nonlinear multidimensional parabolic equations}, Indiana
  Univ. Math. J. {49} (2000), pp.~723--749.

\bibitem{Bowman_1984}
A.~W. Bowman, \emph{An alternative method of cross-validation for the smoothing
  of density estimates}, Biometrika {71} (1984), pp.~353--360.

\bibitem{5}
H.~Brezis and M.~G. Crandall, \emph{Uniqueness of solutions of the
  initial-value problem for {$u\sb{t}-\Delta \varphi (u)=0$}}, J. Math. Pures
  Appl. (9) {58} (1979), pp.~153--163.

\bibitem{2}
R.~Cafiero, V.~Loreto, L.~Pietronero, A.~Vespignani and S.~Zapperi, \emph{Local
  regidity and self-organized criticality for avalanches}, Europhysics Letters
  {29} (1995), pp.~111--116.

\bibitem{Cald_Pulvi}
P.~Calderoni and M.~Pulvirenti, \emph{Propagation of chaos for {B}urgers'
  equation}, Ann. Inst. H. Poincar\'e Sect. A (N.S.) {39} (1983), pp.~85--97.

\bibitem{17}
F.~Cavalli, G.~Naldi, G.~Puppo and M.~Semplice, \emph{High-order relaxation
  schemes for nonlinear degenerate diffusion problems}, SIAM J. Numer. Anal.
  {45} (2007), pp.~2098--2119 (electronic).

\bibitem{TR_Cuvelier}
F.~Cuvelier, \emph{Implementing Kernel Density Estimation on GPU: application
  to a probabilistic algorithm for PDEs of porous media type}, Technical
  report. In preparation.

\bibitem{Dawson}
J.~M. Dawson, \emph{Particle simulation of plasmas}, Rev. Modern Phys. {55}
  (1983), pp.~403--447.

\bibitem{Philipow_Figal}
A.~Figalli and R.~Philipowski, \emph{Convergence to the viscous porous medium
  equation and propagation of chaos}, ALEA Lat. Am. J. Probab. Math. Stat. {4}
  (2008), pp.~185--203.

\bibitem{Butcher}
E.~Hairer, S.~P. N{\o}rsett and G.~Wanner, \emph{Solving ordinary differential
  equations. {I}}, second ed, Springer Series in Computational Mathematics~8.
  Springer-Verlag, Berlin, 1993, Nonstiff problems.

\bibitem{18}
A.~Harten and S.~Osher, \emph{Uniformly high-order accurate nonoscillatory
  schemes. {I}}, SIAM J. Numer. Anal. {24} (1987), pp.~279--309.

\bibitem{Hockney_Livre}
R.~W. Hockney and J.~W. Eastwood, \emph{Computer simulation using particles}.
  McGraw-Hill, New York, 1981.

\bibitem{Jin_Levermore}
S.~Jin and C.~D. Levermore, \emph{Numerical schemes for hyperbolic conservation
  laws with stiff relaxation terms}, J. Comput. Phys. {126} (1996),
  pp.~449--467.

\bibitem{Jin}
S.~Jin and Z.~P. Xin, \emph{The relaxation schemes for systems of conservation
  laws in arbitrary space dimensions}, Comm. Pure Appl. Math. {48} (1995),
  pp.~235--276.

\bibitem{Jones_al_1996}
M.~C. Jones, J.~S. Marron and S.~J. Sheather, \emph{A brief survey of bandwidth
  selection for density estimation}, J. Amer. Statist. Assoc. {91} (1996),
  pp.~401--407.

\bibitem{9}
B.~Jourdain, \emph{Probabilistic approximation for a porous medium equation},
  Stochastic Process. Appl. {89} (2000), pp.~81--99.

\bibitem{11}
B.~Jourdain and S.~M{\'e}l{\'e}ard, \emph{Propagation of chaos and fluctuations
  for a moderate model with smooth initial data}, Ann. Inst. H. Poincar\'e
  Probab. Statist. {34} (1998), pp.~727--766.

\bibitem{Kacur}
J.~Ka{\v{c}}ur, A.~Handlovi{\v{c}}ov{\'a} and M.~Ka{\v{c}}urov{\'a},
  \emph{Solution of nonlinear diffusion problems by linear approximation
  schemes}, SIAM J. Numer. Anal. {30} (1993), pp.~1703--1722.

\bibitem{KARSH}
I.~Karatzas and S.~E. Shreve, \emph{Brownian motion and stochastic calculus},
  second ed, Graduate Texts in Mathematics 113. Springer-Verlag, New York,
  1991.

\bibitem{7}
H.~P.~Jr. McKean, \emph{Propagation of chaos for a class of non-linear
  parabolic equations.}, Stochastic {D}ifferential {E}quations ({L}ecture
  {S}eries in {D}ifferential {E}quations, {S}ession 7, {C}atholic {U}niv.,
  1967), Air Force Office Sci. Res., Arlington, Va., 1967, pp.~41--57.

\bibitem{13}
S.~M{\'e}l{\'e}ard and S.~Roelly-Coppoletta, \emph{A propagation of chaos
  result for a system of particles with moderate interaction}, Stochastic
  Process. Appl. {26} (1987), pp.~317--332.

\bibitem{Oelsch}
K.~Oelschl{\"a}ger, \emph{A law of large numbers for moderately interacting
  diffusion processes}, Z. Wahrsch. Verw. Gebiete {69} (1985), pp.~279--322.

\bibitem{Oelsh87}
\bysame, \emph{A fluctuation theorem for moderately interacting diffusion
  processes}, Probab. Theory Related Fields {74} (1987), pp.~591--616.

\bibitem{Oelsch_Sim}
\bysame, \emph{Simulation of the solution of a viscous porous medium equation
  by a particle method}, SIAM J. Numer. Anal. {40} (2002), pp.~1716--1762
  (electronic).

\bibitem{ParRusso}
L.~Pareschi and G.~Russo, \emph{Implicit-{E}xplicit {R}unge-{K}utta schemes and
  applications to hyperbolic systems with relaxation}, J. Sci. Comput. {25}
  (2005), pp.~129--155.

\bibitem{15}
E.~Parzen, \emph{On estimation of a probability density function and mode},
  Ann. Math. Statist. {33} (1962), pp.~1065--1076.

\bibitem{Philipow}
R.~Philipowski, \emph{Interacting diffusions approximating the porous medium
  equation and propagation of chaos}, Stochastic Process. Appl. {117} (2007),
  pp.~526--538.

\bibitem{Pop_Yong}
I.~S. Pop and W.~Yong, \emph{A numerical approach to degenerate parabolic
  equations}, Numer. Math. {92} (2002), pp.~357--381.

\bibitem{Rudemo_1982}
M.~Rudemo, \emph{Empirical choice of histograms and kernel density estimators},
  Scand. J. Statist. {9} (1982), pp.~65--78.

\bibitem{Scott_Terrell_1987}
D.~W. Scott and G.~R. Terrell, \emph{Biased and unbiased cross-validation in
  density estimation}, J. Amer. Statist. Assoc. {82} (1987), pp.~1131--1146.

\bibitem{SJ}
S.~J. Sheather and M.~C. Jones, \emph{A reliable data-based bandwidth selection
  method for kernel density estimation}, J. Roy. Statist. Soc. Ser. B {53}
  (1991), pp.~683--690.

\bibitem{Show}
R.~E. Showalter, \emph{Monotone operators in {B}anach space and nonlinear
  partial differential equations}, Mathematical Surveys and Monographs~49.
  American Mathematical Society, Providence, RI, 1997.

\bibitem{Shu}
C.~Shu, \emph{Essentially non-oscillatory and weighted essentially
  non-oscillatory schemes for hyperbolic conservation laws}, Advanced numerical
  approximation of nonlinear hyperbolic equations ({C}etraro, 1997), Lecture
  Notes in Math. 1697, Springer, Berlin, 1998, pp.~325--432.

\bibitem{14}
B.~W. Silverman, \emph{Density estimation for statistics and data analysis},
  Monographs on Statistics and Applied Probability. Chapman \& Hall, London,
  1986.

\bibitem{StrVar}
D.~W. Stroock and S.~R.~S. Varadhan, \emph{Multidimensional diffusion
  processes}, Classics in Mathematics. Springer-Verlag, Berlin, 2006, Reprint
  of the 1997 edition.

\bibitem{12}
A.~S. Sznitman, \emph{Topics in propagation of chaos}, \'{E}cole d'\'{E}t\'e de
  {P}robabilit\'es de {S}aint-{F}lour {XIX}---1989, Lecture Notes in Math.
  1464, Springer, Berlin, 1991, pp.~165--251.

\bibitem{Terrel_1990}
G.~R. Terrell, \emph{The maximal smoothing principle in density estimation}, J.
  Amer. Statist. Assoc. {85} (1990), pp.~470--477.

\bibitem{Terrell_Scott_1995}
G.~R. Terrell and D.~W. Scott, \emph{Oversmoothed nonparametric density
  estimates}, J. Amer. Statist. Assoc. {80} (1985), pp.~209--214.

\bibitem{Wand_Jones_1995}
M.~P. Wand and M.~C. Jones, \emph{Kernel smoothing}, Monographs on Statistics
  and Applied Probability~60. Chapman and Hall Ltd., London, 1995.

\bibitem{Woodroofe}
M.~Woodroofe, \emph{On choosing a delta-sequence}, Ann. Math. Statist. {41}
  (1970), pp.~1665--1671.

\bibitem{BotevPreparation}
J.~F.~Grotowski Z.~I.~Botev and D.~P. Kroese, \emph{Kernel density estimation
  via diffusion}, Submitted to the Annals of statistics. (2007).

\end{thebibliography}
